\documentclass[12pt]{article}
\usepackage{amssymb,amsfonts,amsmath, psfrag,eepic,colordvi,graphicx,epsfig,ytableau}
\usepackage[enableskew]{youngtab}
\topskip  
\numberwithin{equation}{section}
\newtheorem{theorem}{Theorem}[section]

\newtheorem{conjecture}[theorem]{Conjecture}

\newtheorem{lemma}[theorem]{Lemma}

\begin{document}
\parskip 6pt

\pagenumbering{arabic}
\def\sof{\hfill\rule{2mm}{2mm}}
\def\ls{\leq}
\def\gs{\geq}
\def\SS{\mathcal S}
\def\qq{{\bold q}}
\def\MM{\mathcal M}
\def\TT{\mathcal T}
\def\EE{\mathcal E}
\def\lsp{\mbox{lsp}}
\def\rsp{\mbox{rsp}}
\def\pf{\noindent {\it Proof.} }
\def\mp{\mbox{pyramid}}
\def\mb{\mbox{block}}
\def\mc{\mbox{cross}}
\def\qed{\hfill \rule{4pt}{7pt}}
\def\block{\hfill \rule{5pt}{5pt}}

\begin{center}
{\Large\bf  On  $(2k+1, 2k+3)$-core partitions with distinct parts}
\end{center}

\begin{center}
{\small Sherry H.F. Yan$^1$, Guizhi Qin$^2$, Zemin Jin$^3$, Robin D.P. Zhou$^4$}

$^{1,2,3}$Department of Mathematics\\
Zhejiang Normal University\\
 Jinhua 321004, P.R. China

$^4$College of Mathematics Physics and Information\\
Shaoxing University\\
Shaoxing 312000, P.R. China

 $^1$huifangyan@hotmail.com

\end{center}

\noindent {\bf Abstract.} In this paper, we are mainly concerned with the enumeration of $(2k+1, 2k+3)$-core partitions with distinct parts. We derive  the number and the largest size of  such partitions, confirming two conjectures  posed by Straub.




\section{Introduction}
 A {\em partition} $\lambda$ of a positive integer $n$ is defined to
be a sequence of nonnegative  integers $(\lambda_1, \lambda_2, \ldots, \lambda_m)$ such that $\lambda_1+\lambda_2+\cdots+
\lambda_m=n$ and $\lambda_1\geq \lambda_2\cdots \geq \lambda_m$.
The empty partition is denoted by $\emptyset$.
 We write $\lambda=(\lambda_1, \lambda_2, \ldots, \lambda_m)\vdash n$ and we say that $n$ is the size of $\lambda$, denoted by $|\lambda|$.
 The {\em Young diagram} of $\lambda$ is defined to be an up- and left- justified array of $n$ boxes with $\lambda_i$ boxes in the $i$-th row. The {\em hook} of each box $B$ in $\lambda$ consists of the box $B$ itself and boxes directly to the right and directly below $B$. The {\em hook length} of $B$, denoted by $h(B)$, is the number of boxes in the hook of $B$.

For a partition $\lambda$, the $\beta$-set of $\lambda$, denoted by $\beta(\lambda)$, is defined to the set of hook lengths of the boxes in the first column of $\lambda$.  For example, Figure \ref{Young} illustrates the Young diagram and the  hook lengths of a partition $\lambda=(5,3,3,1)$.
The $\beta$-set of $\lambda$ is $\beta(\lambda)=\{8,5,4,1\}$.
 Notice  that a partition $\lambda$ is uniquely determined by its  $\beta$-set. Given a decreasing sequence of positive integers  $(h_1, h_2, \ldots, h_m)$, it is easily seen that the unique partition $\lambda$ with $\beta(\lambda)=\{h_1, h_2, \ldots, h_m\}$ is
$$
\lambda=(h_1-(m-1), h_2-(m-2), \ldots, h_{m-1}-1, h_m).
$$
Hence, we have

\begin{equation}\label{eq}
|\lambda|=\sum_{i=1}^{m}h_m-{m\choose 2}.
\end{equation}

\begin{figure}[h]
\begin{center}
 \begin{ytableau}
     8& 6 & 5 & 2&1\\
     5 & 3 & 2\\
    4 & 2 &1\\
     1
 \end{ytableau}
 \end{center}
\caption{ The Young  diagram  of $\lambda=(5,3,3,1)$.}\label{Young}
\end{figure}

For a positive integer $t$, a partition is said to be a {\em $t$-core partition}, or simply a {\em $t$-core},  if it contains no box whose hook length is a multiple of $t$. Let $s$ be a positive integer not equal to $t$, we say that $\lambda$ is an  {\em $(s,t)$-core partition} if it is simultaneously an  $s$-core and a $t$-core. Anderson \cite{And}  showed that the number of $(s,t)$-core partitions is the
rational Catalan number  ${1\over s+t}{s+t\choose s}$
when $s$ and $t$ are coprime to each other.
The proof of Anderson's theorem is through  characterizing the $\beta$-sets of $(s,t)$-core partitions as
order ideals of the poset $P_{s,t}$,  where
\[P_{(s,t)}=\mathbb{N}^+\setminus \{n \in\mathbb{N}^+\mid n=k_1s+k_2t \mbox{ for some } k_1,k_2\in \mathbb{N}\}\]
  whose partial order is fixed by requiring $x\in P_{s,t}$ to cover
$y\in P_{s,t}$ if $x-y$ is either $s$ or $t$.

  Simultaneous core partitions have been   extensively studied. Results on the number, the largest size and the
average size of such partitions could be found in \cite{Agg, Amd-lev, Arm, Chen, Joh, Sta2, Wang, Xiong2, Yang}.
   Recently, Straub \cite{Str} and  Xiong \cite{Xiong} independently  showed that the number of $(s,s+1)$-core  partitions into distinct parts is given by the Fibonacci number $F_{s+1}$, which verify a conjecture posed by Amdeberhan \cite{Amd}. In \cite{Xiong}, Xiong also obtained results on the  the largest size and the average size of such
partitions, which completely settle  Amdeberhan's  conjecture concerning the enumeration of $(s,s+1)$-core  partitions into distinct parts.

In this paper, we are mainly concerned with the number and the largest size of $(2k+1, 2k+3)$-core partition into distinct parts, which verify the following two conjectures posed by
Straub \cite{Str}.

\begin{conjecture}\label{con1}
If $s$ is odd, then the number of $(s, s+2)$-core partitions into distinct parts equals $2^{s-1}$.
\end{conjecture}

\begin{conjecture}\label{con2}
If $s$ is odd, then the largest  size of  $(s,s+2)$-core partitions into distinct parts is given by
${(s^2-1)(s+3)(5s+17)\over 384}$.
\end{conjecture}

\section{Proof of Conjecture \ref{con1}}

In this section,  we aim to confirm Conjecture \ref{con1}. We begin with some definitions and notations.

Throughout the article, we will follow the poset terminology given
by Stanley \cite{Sta}.
Let $P$ be a poset.
For two elements $x$ and $y$ in $P$, we say that $y$ {\em covers} $x$ if $x\prec y$ and there exists no element $z\in P$ such that $x\prec z\prec y$.
Let $P$ be a graded poset.
An element $x$ in $P$ is said to be of {\em rank} $s$ if it  covers an element of rank $s-1$. Note that the elements of rank $0$   in $P$ are just the minimal elements.
The {\em Hasse diagram} of a finite poset $P$ is a graph whose vertices are the elements of $P$,  whose edges are the cover relations, and such that if $y$ covers $x$ then there is an edge connecting $x$ and $y$ and
 $y$ is placed above $x$. For example, in the Hasse diagram of $P_{t,t+1}$,
 each element $x$ of rank $s $ covers exactly two elements $x-t$ and $x-t-1$ of rank $s-1$ for all $1\leq s\leq t-1$, and the elements of rank $0$ are $1,2,\ldots, t-1$.
 See Figure \ref{P{7,8}} for an illustration of
   the Hasse diagram of the poset $P_{7,8}$.  An order ideal of $P$ is a subset $I$ such that if any $y\in I$ and $x\prec y$ in $P$, then $x\in I$.
 Denote by $J(P)$ the set of all order ideals of $P$.

 \begin{figure}[h,t]
\begin{center}
\begin{picture}(60,40)
\setlength{\unitlength}{6mm}

 \put(0,0){\circle*{0.2}}
 \put(0,0){\line(1,1){1}}\put(1,1){\circle*{0.2}}
 \put(1,1){\line(1,-1){1}}\put(2,0){\circle*{0.2}}
 \put(2,0){\line(1,1){1}} \put(2,0){\line(1,1){1}}\put(3,1){\circle*{0.2}}
 \put(3,1){\line(1,-1){1}}\put(4,0){\circle*{0.2}}
 \put(4,0){\line(1,1){1}}\put(5,1){\circle*{0.2}}
 \put(5,1){\line(1,-1){1}}\put(6,0){\circle*{0.2}}
 \put(6,0){\line(1,1){1}}\put(7,1){\circle*{0.2}}
 \put(7,1){\line(1,-1){1}}\put(8,0){\circle*{0.2}}
 \put(8,0){\line(1,1){1}}\put(9,1){\circle*{0.2}}
 \put(9,1){\line(1,-1){1}}\put(10,0){\circle*{0.2}}

 \put(2,2){\circle*{0.2}}\put(2,2){\line(1,-1){1}}
 \put(2,2){\line(-1,-1){1}}\put(4,2){\circle*{0.2}}\put(4,2){\line(1,-1){1}}
 \put(4,2){\line(-1,-1){1}}\put(6,2){\circle*{0.2}}\put(6,2){\line(1,-1){1}}
 \put(6,2){\line(-1,-1){1}}\put(8,2){\circle*{0.2}}\put(8,2){\line(1,-1){1}}
 \put(8,2){\line(-1,-1){1}}

 \put(3,3){\circle*{0.2}}\put(3,3){\line(1,-1){1}}
 \put(3,3){\line(-1,-1){1}}\put(5,3){\circle*{0.2}}\put(5,3){\line(1,-1){1}}
 \put(5,3){\line(-1,-1){1}}\put(7,3){\circle*{0.2}}\put(7,3){\line(1,-1){1}}
 \put(7,3){\line(-1,-1){1}}

 \put(4,4){\circle*{0.2}}\put(4,4){\line(1,-1){1}}
 \put(4,4){\line(-1,-1){1}}\put(6,4){\circle*{0.2}}\put(6,4){\line(1,-1){1}}
 \put(6,4){\line(-1,-1){1}}

 \put(5,5){\circle*{0.2}}\put(5,5){\line(1,-1){1}}
 \put(5,5){\line(-1,-1){1}}
 \put(-0.5, -0.2){\small$1$}\put(1.5, -0.2){\small$2$}\put(3.5, -0.2){\small$3$}\put(5.5, -0.2){\small$4$}\put(7, -0.2){\small$5$}
 \put(9, -0.2){\small$6$}

 \put(0, 0.8){\small$9$}\put(2, 0.8){\small$10$}\put(4, 0.8){\small$11$}\put(6, 0.8){\small$12$}\put(8, 0.8){\small$13$}

  \put(1, 1.8){\small$17$}\put(3, 1.8){\small$18$}\put(5, 1.8){\small$19$}\put(7, 1.8){\small$20$}

  \put(2, 2.8){\small$25$}\put(4, 2.8){\small$26$}\put(6, 2.8){\small$27$}

  \put(3, 3.8){\small$33$}\put(5, 3.8){\small$34$}
  \put(4, 4.8){\small$41$}

\end{picture}
\end{center}
\caption{ The Hasse diagram of the poset  $P_{7,8}$.}\label{P{7,8}}
\end{figure}
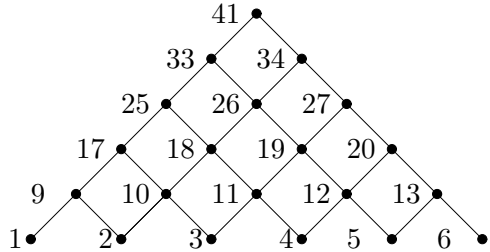

In the following  lemma, we give a characterization of the poset $P_{2k+1, 2k+3}$.

\begin{lemma}\label{mainlem1}
 Let $k$ be a positive integer.
  Then we have \[P_{2k+1, 2k+3}=Q_k\cup Q'_k ,\] where
  \[Q_k=\{2i-1\mid 1\leq i\leq k\}\]
  and \[Q'_k=\{2i+(s-1)(2k+3)\mid 1\leq s\leq 2k, 1\leq i\leq 2k+1-s\}.\]

\end{lemma}

\pf
First we aim to show that  $Q_k\cup Q'_k\subseteq P_{2k+1, 2k+3}$.  Assume to the contrary that  there exists an element $y\in Q_k\cup Q'_k$ such that $y$ is not contained in $P_{2k+1, 2k+3}$.
Choose $y$ to be the smallest such element. Clearly, we have $y=2i+(s-1)(2k+3)$ for some $1\leq s\leq 2k$ and  $1\leq i\leq 2k+1-s$.
 Since $y\not \in P_{2k+1, 2k+3}$,
we have that $2i+(s-1)(2k+3)=m(2k+1)+n(2k+3)$
for some nonnegative integers $m$ and $n$.
We asset that $n=0$.
Otherwise, if $n>0$, then $s>1$.
It follows that $y-(2k+3)=2i+(s-2)(2k+3)\in Q'_k$ and
 $y-(2k+3)=m(2k+1)+(n-1)(2k+3)\not\in P_{2k+1, 2k+3}$.
This contradicts the fact that
$y$ is the smallest element that we choose.
Hence we obtain that $y=2i+(s-1)(2k+3)=m(2k+1)$.
That is $2i+2s-2=(m-s+1)(2k+1)$. Clearly, we have $m>s+2$. Recall that $i\leq 2k+1-s$. This yields that
$2i+2s-2\leq 4k$. Using the fact that $m>s+2$, we have $(m-s+1)(2k+1)\geq 4k+2$. This implies that $2i+2s-2<(m-s+1)(2k+1)$, which yields a contradiction with the relation $2i+2s-2=(m-s+1)(2k+1)$. Hence, we have  $Q_k\cup Q'_k\subseteq P_{2k+1, 2k+3}$.

To complete the proof of the lemma, it remains to
show that $P_{2k+1, 2k+3}\subseteq Q_k\cup Q'_k$.
 Let $y$ be an element in $P_{2k+1, 2k+3}$.
We proceed to show that $y \in Q_k\cup Q'_k$.
It is apparent that $y=x+p(2k+3)$ for some $1\leq x\leq 2k+2$, $x\neq 2k+1$ and $p\geq 0$.
We consider two cases.

\noindent Case 1: $x$ is odd.
That is $y=2i-1+p(2k+3)$ for some $1\leq i\leq k$ and  $p\geq 0$.
If $p=0$, then we have $y\in Q_k$. Otherwise, we have
$y=2i-1+p(2k+3)=2(k+i+1)+(p-1)(2k+3)$.  By setting  $m=k+i+1$, we have
$y=2m+(p-1)(2k+3)$.  In order to show that $y\in Q'_k$, it suffices to verify that
 $ p\leq 2k$ and
   $m\leq 2k+1-p$.
  Assume to the contrary that   $p\geq 2k+2-m$.
Since $P_{2k+1, 2k+3}$ is an order ideal and $y\in P_{2k+1, 2k+3}$,
we have $2m+(2k+1-m)(2k+3)=(2k+3-m)(2k+1)\in P_{2k+1, 2k+3}$.
This contradicts the definition of $P_{2k+1, 2k+3}$. Thus, we have $m\leq 2k+1-p$. Since $m=k+i+1\geq 1$, we have $m\geq 1$. This yields that $p\leq 2k$.

\noindent Case 2: $x$ is even.
That is $y=2i+p(2k+3)  $ for some $1\leq i\leq k+1$ and $p\geq 0$. In order to show that $y\in Q'_k$, it remains to verify that $p\leq 2k-1$ and $i\leq 2k-p$.
Assume to the contrary that  $p\geq 2k-i+1$.
By the definition of the poset $P_{2k+1, 2k+3}$ and
$y\in P_{2k+1, 2k+3}$, we have that
$2i+(2k-i+1)(2k+3)= (2k+3-i)(2k+1)\in P_{2k+1, 2k+3}$. This contradicts the definition of the poset $P_{2k+1, 2k+3}$. Hence we have concluded that $i\leq 2k-p$. Using the fact that $i\geq 1$, we have $p\leq 2k-1$.

Combining the above two cases, we obtain that
$P_{2k+1, 2k+3}\subseteq Q_k\cup Q'_k $.
This completes the proof.\qed

Denote by $A_s=\{2i+(s-1)(2k+3)\mid   1\leq i\leq 2k+1-s\} $  for $1\leq s\leq 2k$.
According to Lemma \ref{mainlem1},  we have $P_{2k+1, 2k+3}=Q_k\cup A_1\cup A_2\cup \ldots \cup A_{2k}$. It is easily seen that  each element $x$ of $A_{s+1}$ covers exactly two elements $x-(2k+3)$ and $x-(2k+1)$ in $A_s$ for all $1\leq s\leq 2k-1$. Moreover, each element $2k+2i+2$ covers exactly two elements $2i-1$ and $2i+1$ for all $1\leq i\leq k-1$, and the element $2k+2$ covers the element $1$.   Hence, the Hasse diagram of the poset $P_{2k+1, 2k+3}$ can be easily constructed by obeying the above rules. For example, the Hasse diagram of the poset $P_{13,15}$ is shown in Figure \ref{P{13,15}}.

\begin{figure}[h,t]
\begin{center}
\begin{picture}(120,80)
\setlength{\unitlength}{6mm}

 \put(0,1){\circle*{0.2}}
 \put(0,1){\line(1,1){1}}

 \put(1,2){\circle*{0.2}}
 \put(1,2){\line(1,-1){1}}
 \put(2,1){\circle*{0.2}}
 \put(2,1){\line(1,1){1}}
 \put(3,2){\circle*{0.2}}
 \put(3,2){\line(1,-1){1}}

 \put(4,1){\circle*{0.2}}
 \put(4,1){\line(1,1){1}}
 \put(5,2){\circle*{0.2}}
 \put(5,2){\line(1,-1){1}}

 \put(6,1){\circle*{0.2}}
 \put(6,1){\line(1,1){1}}\put(7,2){\circle*{0.2}}
 \put(7,2){\line(1,-1){1}}\put(8,1){\circle*{0.2}}
 \put(8,1){\line(1, 1){1}}\put(9,2){\circle*{0.2}}
 \put(9,2){\line(1,-1){1}}\put(10,1){\circle*{0.2}}
 \put(10,1){\line(1,1){1}}\put(11,2){\circle*{0.2}}
 \put(11,2){\line(1,-1){1}}\put(12,1){\circle*{0.2}}
 \put(12,1){\line(1,1){1}}\put(13,2){\circle*{0.2}}
 \put(13,2){\line(1,-1){1}}\put(14,1){\circle*{0.2}}
 \put(14,1){\line(1,1){1}}\put(15,2){\circle*{0.2}}
 \put(15,2){\line(1,-1){1}}\put(16,1){\circle*{0.2}}
 \put(16,1){\line(1, 1){1}}\put(17,2){\circle*{0.2}}
 \put(17,2){\line(1, -1){1}}\put(18,1){\circle*{0.2}}
 \put(18,1){\line(1, 1){1}}\put(19,2){\circle*{0.2}}
 \put(19,2){\line(1, -1){1}}\put(20,1){\circle*{0.2}}
 \put(20,1){\line(1, 1){1}}\put(21,2){\circle*{0.2}}
 \put(21,2){\line(1, -1){1}}\put(22,1){\circle*{0.2}}
   \put(2,3){\circle*{0.2}}\put(2,3){\line(-1, -1){1}}
 \put(2,3){\line(1, -1){1}}

 \put(4,3){\circle*{0.2}}\put(4,3){\line(-1, -1){1}}
 \put(4,3){\line(1, -1){1}}

  \put(6,3){\circle*{0.2}}\put(6,3){\line(-1, -1){1}}
 \put(6,3){\line(1, -1){1}}
 \put(8,3){\circle*{0.2}}\put(8,3){\line(-1, -1){1}}
 \put(8,3){\line(1, -1){1}}
 \put(10,3){\circle*{0.2}}\put(10,3){\line(-1, -1){1}}
 \put(10,3){\line(1, -1){1}}
 \put(12,3){\circle*{0.2}}\put(12,3){\line(-1, -1){1}}
 \put(12,3){\line(1, -1){1}}
 \put(14,3){\circle*{0.2}}\put(14,3){\line(-1, -1){1}}
 \put(14,3){\line(1, -1){1}}
 \put(16,3){\circle*{0.2}}\put(16,3){\line(-1, -1){1}}
 \put(16,3){\line(1, -1){1}}
 \put(18,3){\circle*{0.2}}\put(18,3){\line(-1, -1){1}}
 \put(18,3){\line(1, -1){1}}
 \put(20,3){\circle*{0.2}}\put(20,3){\line(-1, -1){1}}
 \put(20,3){\line(1, -1){1}}

 \put(3,4){\circle*{0.2}}\put(3,4){\line(-1, -1){1}}
 \put(3,4){\line(1, -1){1}}

 \put(5,4){\circle*{0.2}}\put(5,4){\line(-1, -1){1}}
 \put(5,4){\line(1, -1){1}}

  \put(7,4){\circle*{0.2}}\put(7,4){\line(-1, -1){1}}
 \put(7,4){\line(1, -1){1}}
 \put(9,4){\circle*{0.2}}\put(9,4){\line(-1, -1){1}}
 \put(9,4){\line(1, -1){1}}
 \put(11,4){\circle*{0.2}}\put(11,4){\line(-1, -1){1}}
 \put(11,4){\line(1, -1){1}}
 \put(13,4){\circle*{0.2}}\put(13,4){\line(-1, -1){1}}
 \put(13,4){\line(1, -1){1}}
 \put(15,4){\circle*{0.2}}\put(15,4){\line(-1, -1){1}}
 \put(15,4){\line(1, -1){1}}
 \put(17,4){\circle*{0.2}}\put(17,4){\line(-1, -1){1}}
 \put(17,4){\line(1, -1){1}}
 \put(19,4){\circle*{0.2}}\put(19,4){\line(-1, -1){1}}
 \put(19,4){\line(1, -1){1}}

 \put(4,5){\circle*{0.2}}\put(4,5){\line(-1, -1){1}}
 \put(4,5){\line(1, -1){1}}

 \put(6,5){\circle*{0.2}}\put(6,5){\line(-1, -1){1}}
 \put(6,5){\line(1, -1){1}}

  \put(8,5){\circle*{0.2}}\put(8,5){\line(-1, -1){1}}
 \put(8,5){\line(1, -1){1}}
 \put(10,5){\circle*{0.2}}\put(10,5){\line(-1, -1){1}}
 \put(10,5){\line(1, -1){1}}
 \put(12,5){\circle*{0.2}}\put(12,5){\line(-1, -1){1}}
 \put(12,5){\line(1, -1){1}}
 \put(14,5){\circle*{0.2}}\put(14,5){\line(-1, -1){1}}
 \put(14,5){\line(1, -1){1}}
 \put(16,5){\circle*{0.2}}\put(16,5){\line(-1, -1){1}}
 \put(16,5){\line(1, -1){1}}
 \put(18,5){\circle*{0.2}}\put(18,5){\line(-1, -1){1}}
 \put(18,5){\line(1, -1){1}}

  \put(5,6){\circle*{0.2}}\put(5,6){\line(-1, -1){1}}
 \put(5,6){\line(1, -1){1}}

 \put(7,6){\circle*{0.2}}\put(7,6){\line(-1, -1){1}}
 \put(7,6){\line(1, -1){1}}

  \put(9,6){\circle*{0.2}}\put(9,6){\line(-1, -1){1}}
 \put(9,6){\line(1, -1){1}}
 \put(11,6){\circle*{0.2}}\put(11,6){\line(-1, -1){1}}
 \put(11,6){\line(1, -1){1}}
 \put(13,6){\circle*{0.2}}\put(13,6){\line(-1, -1){1}}
 \put(13,6){\line(1, -1){1}}
 \put(15,6){\circle*{0.2}}\put(15,6){\line(-1, -1){1}}
 \put(15,6){\line(1, -1){1}}
 \put(17,6){\circle*{0.2}}\put(17,6){\line(-1, -1){1}}
 \put(17,6){\line(1, -1){1}}

 \put(6,7){\circle*{0.2}}\put(6,7){\line(-1, -1){1}}
 \put(6,7){\line(1, -1){1}}

 \put(8,7){\circle*{0.2}}\put(8,7){\line(-1, -1){1}}
 \put(8,7){\line(1, -1){1}}

  \put(10,7){\circle*{0.2}}\put(10,7){\line(-1, -1){1}}
 \put(10,7){\line(1, -1){1}}
 \put(12,7){\circle*{0.2}}\put(12,7){\line(-1, -1){1}}
 \put(12,7){\line(1, -1){1}}
 \put(14,7){\circle*{0.2}}\put(14,7){\line(-1, -1){1}}
 \put(14,7){\line(1, -1){1}}
 \put(16,7){\circle*{0.2}}\put(16,7){\line(-1, -1){1}}
 \put(16,7){\line(1, -1){1}}

 \put(7,8){\circle*{0.2}}\put(7,8){\line(-1, -1){1}}
 \put(7,8){\line(1, -1){1}}

 \put(9,8){\circle*{0.2}}\put(9,8){\line(-1, -1){1}}
 \put(9,8){\line(1, -1){1}}

  \put(11,8){\circle*{0.2}}\put(11,8){\line(-1, -1){1}}
 \put(11,8){\line(1, -1){1}}
 \put(13,8){\circle*{0.2}}\put(13,8){\line(-1, -1){1}}
 \put(13,8){\line(1, -1){1}}
 \put(15,8){\circle*{0.2}}\put(15,8){\line(-1, -1){1}}
 \put(15,8){\line(1, -1){1}}

 \put(8,9){\circle*{0.2}}\put(8,9){\line(-1, -1){1}}
 \put(8,9){\line(1, -1){1}}

 \put(10,9){\circle*{0.2}}\put(10,9){\line(-1, -1){1}}
 \put(10,9){\line(1, -1){1}}

  \put(12,9){\circle*{0.2}}\put(12,9){\line(-1, -1){1}}
 \put(12,9){\line(1, -1){1}}
 \put(14,9){\circle*{0.2}}\put(14,9){\line(-1, -1){1}}
 \put(14,9){\line(1, -1){1}}

  \put(9,10){\circle*{0.2}}\put(9,10){\line(-1, -1){1}}
 \put(9,10){\line(1, -1){1}}

 \put(11,10){\circle*{0.2}}\put(11,10){\line(-1, -1){1}}
 \put(11,10){\line(1, -1){1}}

  \put(13,10){\circle*{0.2}}\put(13,10){\line(-1, -1){1}}
 \put(13,10){\line(1, -1){1}}

 \put(10,11){\circle*{0.2}}\put(10,11){\line(-1, -1){1}}
 \put(10,11){\line(1, -1){1}}

 \put(12,11){\circle*{0.2}}\put(12,11){\line(-1, -1){1}}
 \put(12,11){\line(1, -1){1}}

  \put(11,12){\circle*{0.2}}\put(11,12){\line(-1, -1){1}}
 \put(11,12){\line(1, -1){1}}

  \put(12,1){\line(1, -1){1}}\put(13,0){\circle*{0.2}}\put(13,0){\line(1, 1){1}}
  \put(15,0){\line(-1, 1){1}}\put(15,0){\circle*{0.2}}\put(15,0){\line(1, 1){1}}\put(17,0){\line(-1, 1){1}}\put(17,0){\circle*{0.2}}\put(17,0){\line(1, 1){1}}\put(19,0){\line(-1, 1){1}}\put(19,0){\circle*{0.2}}\put(19,0){\line(1, 1){1}}\put(21,0){\line(-1, 1){1}}\put(21,0){\circle*{0.2}}\put(21,0){\line(1, 1){1}} \put(23,0){\line(-1, 1){1}}\put(23,0){\circle*{0.2}}

 \put(-0.6, 0.8){\small$2$}\put(1.4, 0.8){\small$4$}\put(3.4, 0.8){\small$6$}
 \put(5.4, 0.8){\small$8$} \put(7, 0.8){\small$10$}\put(9, 0.8){\small$12$}
 \put(11, 0.8){\small$14$}\put(13, 0.8){\small$16$}\put(15, 0.8){\small$18$}
 \put(17, 0.8){\small$20$}\put(19, 0.8){\small$22$} \put(21, 0.8){\small$24$}

 \put(0, 1.8){\small$17$}\put(2, 1.8){\small$19$}\put(4, 1.8){\small$21$}
 \put(6, 1.8){\small$23$} \put(8, 1.8){\small$25$}\put(10, 1.8){\small$27$}
 \put(12, 1.8){\small$29$}\put(14, 1.8){\small$31$}\put(16, 1.8){\small$33$}
 \put(18, 1.8){\small$35$}\put(20, 1.8){\small$37$}

 \put(1, 2.8){\small$32$}\put(3, 2.8){\small$34$}\put(5, 2.8){\small$36$}
 \put(7, 2.8){\small$38$} \put(9, 2.8){\small$40$}\put(11, 2.8){\small$42$}
 \put(13, 2.8){\small$44$}\put(15, 2.8){\small$46$}\put(17, 2.8){\small$48$}
 \put(19, 2.8){\small$50$}

 \put(2, 3.8){\small$47$}\put(4, 3.8){\small$49$}\put(6, 3.8){\small$51$}
 \put(8, 3.8){\small$53$} \put(10, 3.8){\small$55$}\put(12, 3.8){\small$57$}
 \put(14, 3.8){\small$59$}\put(16, 3.8){\small$61$}\put(18, 3.8){\small$63$}

  \put(3, 4.8){\small$62$}\put(5, 4.8){\small$64$}\put(7, 4.8){\small$66$}
 \put(9, 4.8){\small$68$} \put(11, 4.8){\small$70$}\put(13, 4.8){\small$72$}
 \put(15, 4.8){\small$74$}\put(17, 4.8){\small$76$}

 \put(4, 5.8){\small$77$}\put(6, 5.8){\small$79$}\put(8, 5.8){\small$81$}
 \put(10, 5.8){\small$83$} \put(12, 5.8){\small$85$}\put(14, 5.8){\small$87$}
 \put(16, 5.8){\small$89$}

 \put(5, 6.8){\small$92$}\put(6.8, 6.8){\small$94$}\put(8.8, 6.8){\small$96$}
 \put(10.8, 6.8){\small$98$} \put(12.8, 6.8){\small$100$}\put(14.8, 6.8){\small$102$}

  \put(5.8, 7.8){\small$107$}\put(7.8, 7.8){\small$109$}\put(9.8, 7.8){\small$111$}
 \put(11.8, 7.8){\small$113$} \put(13.8, 7.8){\small$115$}

 \put(6.8, 8.8){\small$122$}\put(8.8, 8.8){\small$124$}\put(10.8, 8.8){\small$126$}
 \put(12.8, 8.8){\small$128$}

 \put(7.8, 9.8){\small$137$}\put(9.8, 9.8){\small$139$}\put(11.8, 9.8){\small$141$}

  \put(8.8, 10.8){\small$152$}\put(10.8, 10.8){\small$154$}
  \put(9.8, 11.8){\small$167$}

  \put(12.6, -0.2){\small$1$}\put(14.4, -0.2){\small$3$}\put(16.4, -0.2){\small$5$}\put(18.4, -0.2){\small$7$}\put(20.4, -0.2){\small$9$}
  \put(22.2, -0.2){\small$11$}
 \end{picture}
\end{center}
\caption{ The Hasse diagram of the poset  $P_{13,15}$.}\label{P{13,15}}
\end{figure}
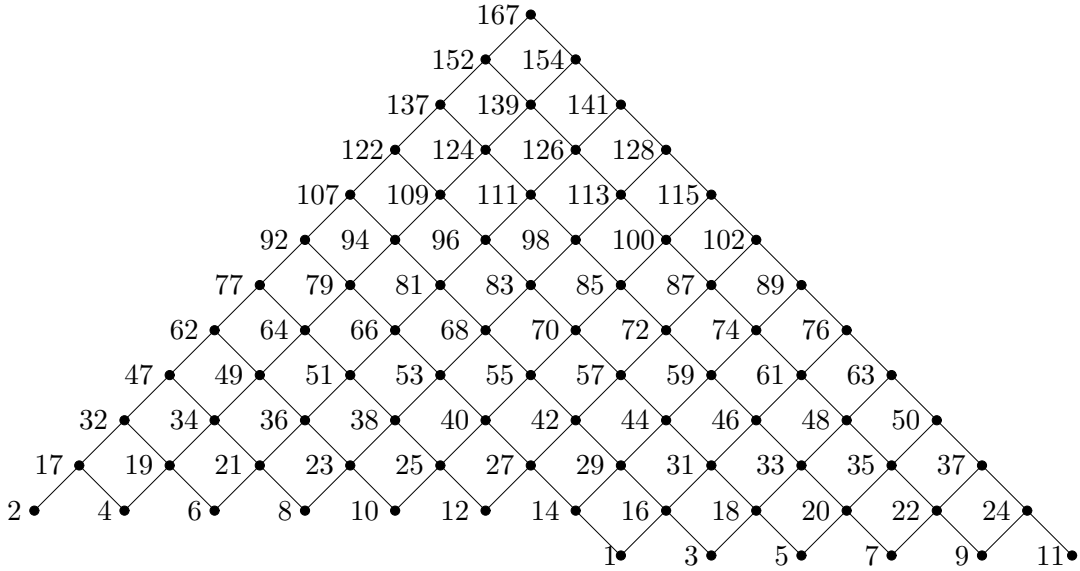

Denote by $M_{2k+1, 2k+3}$   the poset   obtained from  $P_{2k+1, 2k+3}$  by removing  all the elements $y$ such that $y\succeq 2k+2$  and preserving the cover relation among the remaining elements.
 From the Hasse diagram of $P_{2k+1, 2k+3}$,  we immediately obtain the following        characterization of the
poset $M_{2k+1, 2k+3}$.

\begin{lemma}\label{mainlem2}
 Let $k$ be a positive integer.
  Define
  \[T_k=\{2i-1+(s-1)(2k+3)\mid 1\leq s\leq k, 1\leq i\leq k+1-s\}\]
  and
  \[T'_k=\{2i+(s-1)(2k+3)\mid 1\leq s\leq k, 1\leq i\leq k+1-s\}.\]
  Then we have $M_{2k+1, 2k+3}=T_k\cup T'_k$.
\end{lemma}

Relying on Lemma \ref{mainlem2}, we can get  the Hasse diagrams of the poset  $M_{2k+1, 2k+3}$.
 Figure  \ref{M{13,15}} illustrates the Hasse diagram of the poset $M_{13,15}$.

 \begin{figure}[h,t]
\begin{center}
\begin{picture}(120,40)
\setlength{\unitlength}{6mm}

 \put(0,0){\circle*{0.2}}
 \put(0,0){\line(1,1){1}}\put(1,1){\circle*{0.2}}
 \put(1,1){\line(1,-1){1}}\put(2,0){\circle*{0.2}}
 \put(2,0){\line(1,1){1}} \put(2,0){\line(1,1){1}}\put(3,1){\circle*{0.2}}
 \put(3,1){\line(1,-1){1}}\put(4,0){\circle*{0.2}}
 \put(4,0){\line(1,1){1}}\put(5,1){\circle*{0.2}}
 \put(5,1){\line(1,-1){1}}\put(6,0){\circle*{0.2}}
 \put(6,0){\line(1,1){1}}\put(7,1){\circle*{0.2}}
 \put(7,1){\line(1,-1){1}}\put(8,0){\circle*{0.2}}
 \put(8,0){\line(1,1){1}}\put(9,1){\circle*{0.2}}
 \put(9,1){\line(1,-1){1}}\put(10,0){\circle*{0.2}}

 \put(2,2){\circle*{0.2}}\put(2,2){\line(1,-1){1}}
 \put(2,2){\line(-1,-1){1}}\put(4,2){\circle*{0.2}}\put(4,2){\line(1,-1){1}}
 \put(4,2){\line(-1,-1){1}}\put(6,2){\circle*{0.2}}\put(6,2){\line(1,-1){1}}
 \put(6,2){\line(-1,-1){1}}\put(8,2){\circle*{0.2}}\put(8,2){\line(1,-1){1}}
 \put(8,2){\line(-1,-1){1}}

 \put(3,3){\circle*{0.2}}\put(3,3){\line(1,-1){1}}
 \put(3,3){\line(-1,-1){1}}\put(5,3){\circle*{0.2}}\put(5,3){\line(1,-1){1}}
 \put(5,3){\line(-1,-1){1}}\put(7,3){\circle*{0.2}}\put(7,3){\line(1,-1){1}}
 \put(7,3){\line(-1,-1){1}}

 \put(4,4){\circle*{0.2}}\put(4,4){\line(1,-1){1}}
 \put(4,4){\line(-1,-1){1}}\put(6,4){\circle*{0.2}}\put(6,4){\line(1,-1){1}}
 \put(6,4){\line(-1,-1){1}}

 \put(5,5){\circle*{0.2}}\put(5,5){\line(1,-1){1}}
 \put(5,5){\line(-1,-1){1}}
 \put(-0.5, -0.2){\small$2$}\put(1.5, -0.2){\small$4$}\put(3.5, -0.2){\small$6$}\put(5.5, -0.2){\small$8$}\put(7, -0.2){\small$10$}
 \put(9, -0.2){\small$12$}

 \put(0, 0.8){\small$17$}\put(2, 0.8){\small$19$}\put(4, 0.8){\small$21$}\put(6, 0.8){\small$23$}\put(8, 0.8){\small$25$}

  \put(1, 1.8){\small$32$}\put(3, 1.8){\small$34$}\put(5, 1.8){\small$36$}\put(7, 1.8){\small$38$}

  \put(2, 2.8){\small$47$}\put(4, 2.8){\small$49$}\put(6, 2.8){\small$51$}

  \put(3, 3.8){\small$62$}\put(5, 3.8){\small$64$}
  \put(4, 4.8){\small$77$}


  \put(12,0){\circle*{0.2}}
 \put(12,0){\line(1,1){1}}\put(13,1){\circle*{0.2}}
 \put(13,1){\line(1,-1){1}}\put(14,0){\circle*{0.2}}
 \put(14,0){\line(1,1){1}} \put(14,0){\line(1,1){1}}\put(15,1){\circle*{0.2}}
 \put(15,1){\line(1,-1){1}}\put(16,0){\circle*{0.2}}
 \put(16,0){\line(1,1){1}}\put(17,1){\circle*{0.2}}
 \put(17,1){\line(1,-1){1}}\put(18,0){\circle*{0.2}}
 \put(18,0){\line(1,1){1}}\put(19,1){\circle*{0.2}}
 \put(19,1){\line(1,-1){1}}\put(20,0){\circle*{0.2}}
 \put(20,0){\line(1,1){1}}\put(21,1){\circle*{0.2}}
 \put(21,1){\line(1,-1){1}}\put(22,0){\circle*{0.2}}

 \put(14,2){\circle*{0.2}}\put(14,2){\line(1,-1){1}}
 \put(14,2){\line(-1,-1){1}}\put(16,2){\circle*{0.2}}\put(16,2){\line(1,-1){1}}
 \put(16,2){\line(-1,-1){1}}\put(18,2){\circle*{0.2}}\put(18,2){\line(1,-1){1}}
 \put(18,2){\line(-1,-1){1}}\put(20,2){\circle*{0.2}}\put(20,2){\line(1,-1){1}}
 \put(20,2){\line(-1,-1){1}}

 \put(15,3){\circle*{0.2}}\put(15,3){\line(1,-1){1}}
 \put(15,3){\line(-1,-1){1}}\put(17,3){\circle*{0.2}}\put(17,3){\line(1,-1){1}}
 \put(17,3){\line(-1,-1){1}}\put(19,3){\circle*{0.2}}\put(19,3){\line(1,-1){1}}
 \put(19,3){\line(-1,-1){1}}

 \put(16,4){\circle*{0.2}}\put(16,4){\line(1,-1){1}}
 \put(16,4){\line(-1,-1){1}}\put(18,4){\circle*{0.2}}\put(18,4){\line(1,-1){1}}
 \put(18,4){\line(-1,-1){1}}

 \put(17,5){\circle*{0.2}}\put(17,5){\line(1,-1){1}}
 \put(17,5){\line(-1,-1){1}}
 \put(11.5, -0.2){\small$1$}\put(13.5, -0.2){\small$3$}\put(15.5, -0.2){\small$5$}\put(17.5, -0.2){\small$7$}\put(19, -0.2){\small$9$}
 \put(21, -0.2){\small$11$}

 \put(12, 0.8){\small$16$}\put(14, 0.8){\small$18$}\put(16, 0.8){\small$20$}\put(18, 0.8){\small$22$}\put(20, 0.8){\small$24$}

  \put(13, 1.8){\small$31$}\put(15, 1.8){\small$33$}\put(17, 1.8){\small$35$}\put(19, 1.8){\small$37$}

  \put(14, 2.8){\small$46$}\put(16, 2.8){\small$48$}\put(18, 2.8){\small$50$}

  \put(15, 3.8){\small$61$}\put(17, 3.8){\small$63$}
  \put(16, 4.8){\small$76$}
\end{picture}
\end{center}
\caption{ The Hasse diagram of the poset $M_{13,15}$.}\label{M{13,15}}
\end{figure}
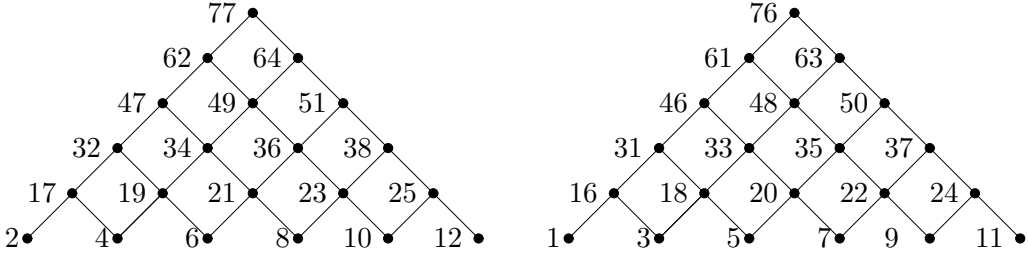

In the following theorem, Anderson \cite{And} established a correspondence between core partitions and order ideals of a certain poset by mapping a partition to its $\beta$-set.
\begin{theorem}\label{th1.1}
Let $s$, $t$ be two coprime positive integers, and let $\lambda$ be a partition of $n$. Then $\lambda$ is an $s$-core (or $(s,t)$-core) partition if and only if
$\beta(\lambda)$ is an order ideal of $P_s$ (or $P_{s,t}$).
\end{theorem}

The following theorem provides a
characterization of the $\beta$-set of partitions into distinct parts \cite{Xiong}.    \begin{theorem}\label{th1.2}
The partition $\lambda$ is a partition into distinct parts if and only if there
doesn't exist $x,y\in \beta(\lambda)$ with $x-y=1$.
   \end{theorem}

  An order ideal $I$ of $P$ is said to be {\em nice} if   there
doesn't exist $x,y\in I$ with $x-y=1$. Denote by $L(P)$ the set of nice order ideals of the poset $P$.   Combining Theorems \ref{th1.1}
and \ref{th1.2}, we obtain the following result.
\begin{theorem}\label{th1.3}
Let $s$, $t$ be two coprime positive integers, and let $\lambda$ be a partition of $n$. Then $\lambda$ is an $s$-core (or $(s,t)$-core) partition into distinct parts if and only if
$\beta(\lambda)$ is a nice  order ideal of $P_s$ (or $P_{s,t}$).
\end{theorem}

In order to get the enumeration of $(2k+1, 2k+3)$-core partitions into distinct parts, we need the notions  of   Dyck   path and free Dyck path.
Recall that a {\em Dyck } path of order $n$ is a lattice path in $\emph{Z}\times \emph{Z}$ from $(0,0)$ to $(2n,0)$
 using up steps $U=(1,1)$ and down steps $D=(1,-1) $, and never
  lying below the $x$-axis \cite{De}.
   Denote by $\mathcal{D}_n$ the set of all  Dyck  paths of order $n$.
   It is well known that Dyck  paths of order $n$ are counted by the $n$-th
    Catalan  number $c_n={1\over n+1}{2n\choose n} $.
     A {\em free } Dyck path is just a Dyck path but
     without the restriction that it
cannot go below the x-axis.
Let $\mathcal{FD}_n$ denote the set of  free Dyck paths
 from $(0,0)$ to $(2n,0)$. By simple arguments we have
  that $|\mathcal{FD}_n|={2n\choose n}$.

We first describe a bijection $\phi$ between the set of order ideals of $P_{t,t+1}$ and the set of Dyck paths from $(0,0)$ to $(2t,0)$. To this end, we shall partition $J(P_{t,t+1})$ according to the smallest missing element  of rank $0$ in an order ideal.   For $1\leq i\leq t-1$, let $J_i(P_{t,t+1})$ denote the set of order ideals of $P_{t,t+1}$ such that $i$ is the smallest missing elements of rank $0$.
Let $J_t(P_{t,t+1})$ denote the set of  order ideals containing all the minimal elements of $P_{t,t+1}$.
Now we are in the position to construct a Dyck path from $I\in J_t(P_{t,t+1})$ recursively as follows.
\begin{itemize}
\item For $2\leq i\leq t-1$ and an  order ideal $I\in J_i(P_{t,t+1})$, we can decompose $I$ into three parts: one is $\{1,2,\ldots, i-1\}$, one is isomorphic to an order ideal $I_1$ of  $J(P_{i-1, i})$, and one is isomorphic to an order ideal $I_2$ of  $J(P_{t-i, t-i+1})$. Define $\phi(I)=U\phi(I_1)D\phi(I_2)$.
    \item  For any order ideal $I\in J_1(P_{t, t+1})$,  $I$ is  isomorphic to an order ideal $I_1$ of $J(P_{t-1, t})$. Define $\phi(I)=UD\phi(I_1)$.
         \item For any order ideal $I\in J_t(P_{t, t+1})$, one can decompose $I$ into two parts: one is ${1,2, \ldots, t-1}$ and  the other is  isomorphic to an order ideal $I_1$ of $J(P_{t-1, t})$. Define $\phi(I)=U \phi(I_1)D$.
\end{itemize}

It is easy to check that the resulting path is a Dyck from $(0,0)$ to $(2t,0)$. Conversely, given a Dyck path $P$, by the first return decomposition, $P$ can be uniquely decomposed as
$$
P=UP_1DP_2,
$$
where $P_1$ and $P_2$ are (possibly empty) Dyck paths. Hence, the above procedure is reversible and the map $\phi$ is a bijection.

 For a path  $P=p_1p_2\ldots p_n$ where $p_i\in \{U,D\}$ for all $1\leq i\leq n$,  the {\em reverse} of the path,
denoted by $\overline{P}$, is defined  by $p_np_{n-1}\ldots p_{1}$.
For example, the reverse of the path  $P=UDDUD$ is given by  $DUDDU$. It is easily seen that for a Dyck path $P$, its reverse $\overline{P}$ is  a free Dyck path which lies weakly below the $x$-axis.
Denote by $\widehat{L}(M_{2k+1, 2k+3})$ the set of order ideals in $L(M_{2k+1, 2k+3})$   which do not contain the element $1$.

\begin{lemma}\label{lem1.1}
Let $k$ be  a nonnegative  integer. There is a bijection between the set   $\widehat{L}(M_{2k+1, 2k+3})$      and the set $\mathcal{FD}_{k}$.
\end{lemma}
\pf First we describe a recursive map $\psi$ from the set $\widehat{L}(M_{2k+1, 2k+3})$    to the set  $\mathcal{FD}_{k}$.  Let $I\in \widehat{L}(M_{2k+1, 2k+3}) $. For $k=0$, let $\psi(I)$ be the empty path.  For $k=1$,   let $\psi(I)=UD$ if $I$ contains the element $2$, and  let $\psi(I)=DU$, otherwise. Assume that for all $2\leq j\leq k-1$ and an order ideal $I'\in \widehat{L}(M_{2j+1, 2j+3})$ which contains the element  $2j$,  we have $\psi(I')\in \mathcal{FD}_{j}$ which ends  with a down step.  Similarly, for  all  $2\leq j\leq k-1$ and an order ideal $I'\in \widehat{L}(M_{2j+1, 2j+3})$ which  does not contain the element $2j$,  we have $\psi(I')\in \mathcal{FD}_{j}$ which ends with an up  step. Now we construct a lattice path from $(0,0)$ to $(2k,0)$ as follows.
\begin{itemize}
 \item  If $I$ contains the element  $2k$, then let $2\ell$ be the largest  missing  even  number  satisfying that $0\leq \ell\leq k-1$.
    In this case, one can decompose  $I$ into three parts: one is $\{2\ell+2, 2\ell+4, \ldots, 2k\}$, one  is isomorphic to an order ideal $I'$ of $J(P_{k-\ell,k-\ell+1})$, and one is  isomorphic to an order ideal $I''$ of $\widehat{L}(M_{2\ell+1, 2\ell+3})$.    It is easily seen that $2\ell\notin I''$.
        Set $\psi(I)=\psi(I'')\phi(I')$.

    \item  If $I$ does not contain the element $2k$  but contains the element $2j$ for some $1\leq j\leq k-1$, then let $\ell$ be the largest  integer such that $I$ contains  the element $2\ell$.
        In this case, one can decompose  $I$ into two parts:   one  is isomorphic to an order ideal $I'$ of $J(P_{k-\ell,k-\ell+1})$, and the other  is isomorphic to an order ideal $I''$ of $\widehat{L}(M_{2\ell+1, 2\ell+3})$.  Clearly, $I''$ contains the element $2\ell$.
         Set $\psi(I)=\psi(I'')\overline{\phi(I')}$.

        \item   If $I$  does not  contain the element  $2j$ for all $1\leq j\leq k$, then $I$  is isomorphic to an order ideal $I'$ of $J(P_{k, k+1})$.  Set $\psi(I)=\overline{\phi(I')}$.
\end{itemize}
By induction hypothesise, one can easily check that the resulting path $\psi(I)$ is a free Dyck path from $(0,0)$ to $(2k,0)$. Moreover, the path $\psi(I)$ ends  with a  down step if $I$ contains the element  $2k$,  and $\psi(I)$ ends with an up  step, otherwise. Thus, the map $\psi$ is well defined, that is, we have $\psi(I)\in \mathcal{FD}_k$ for any $I\in M(P_{2k+1, 2k+3})$.

Conversely, for any free Dyck path $P$ ending with a down  step, it can be uniquely decomposed as
$$
P=P'P''
$$
where  $P' $ is  a (possibly empty) free Dyck path ending with an up step and $P''$  is a Dyck path.
Analogously, for any free Dyck path $P$  ending  with an up step, it can be uniquely decomposed as
$$
P=P'P''
$$
where  $P' $ is  a (possibly empty) free Dyck path ending with a down step and $P''$  is the reverse of a Dyck path.
Hence for any free Dyck path $P\in \mathcal{FD}_k$, we can recover an order ideal in $\widehat{L}(M_{2k+1, 2k+3})$ by reversing the above procedure. This implies that the map $\psi$ is a bijection, which completes the proof.
\qed

By Lemma \ref{lem1.1}, we immediately get the following result on the cardinality  of  the set  $\widehat{L}(M_{2k+1, 2k+3})$.

\begin{lemma}\label{lem1.2}
Let $k$ be a nonnegative  integer. The cardinality of the set  $\widehat{L}(M_{2k+1, 2k+3})$   is given by ${2k\choose k}$.
\end{lemma}

 Let $k\geq 1$. Denote by $\widetilde{L}(M_{2k+1, 2k+3})$ the set of order ideals of $ L(M_{2k+1, 2k+3})$ which contain neither $1$ nor $2k$.
  From the construction of the map $\psi$, it is easy to check that the map $\psi$ sends an order ideal   $I\in \widetilde{L}(M_{2k+1, 2k+3})$  to a free  Dyck path from $(0,0)$ to $(2k,0)$ ending with an up step. The latter is counted by ${2k-1\choose k}$.  Hence, we have the following result.
\begin{lemma}\label{lem1.3}
Let $k$ be a positive  integer. The number of nice order ideals of $\widetilde{L}(M_{2k+1, 2k+3}) $  is counted by ${2k-1\choose k}$.
\end{lemma}

In the following theorem, we shall obtain the enumeration of the cardinality of the set $L(M_{2k+1, 2k+3})$.

\begin{theorem}\label{th1.4}
Let $k$ be a nonnegative integer. The cardinality of the set $L(M_{2k+1, 2k+3})$  is given by $\sum_{i=0}^k{1\over i+1}{2i\choose i} {2k-2i\choose k-i}$.
\end{theorem}
\pf  It is easily seen that the assertion holds for $k=0$.  For $k>0$  and  $0\leq i\leq k$, denote by $L_i(M_{2k+1, 2k+3})$ the set of order  ideals  $I$ of $L(M_{2k+1, 2k+3})$  such that  $2i+1$ is the smallest odd number  missing from $I$. For an order ideal $I\in L_i(M_{2k+1, 2k+3})$ where $1\leq i\leq k$, we can decompose it into three parts: one is  \{1, 3, \ldots, 2i-1\}, one is  isomorphic to an order ideal of $J(P_{i,i+1})$,  and one is  isomorphic to an order ideal of $\widehat{L}(M_{2k-2i+1, 2k-2i+3})$.
Thus, by Lemma \ref{lem1.2},  we deduce that
$$
|L_i(M_{2k+1, 2k+3})|=c_i{2k-2i\choose k-i}={1\over i+1}{2i\choose i} {2k-2i\choose k-i}.
$$

For an order ideal $I\in L_0(M_{2k+1, 2k+3})$, we can decompose it into two parts: one is  \{2, 4, \ldots, 2k\} and the other is isomorphic to an order ideal of $J(P_{k,k+1})$. Thus, it follows that
$$
|L_0(M_{2k+1, 2k+3})|=c_k={1\over k+1}{2k\choose k}.
$$
Hence, we have
$$
|L(M_{2k+1, 2k+3})|=\sum_{i=1}^{k}|L_i(M_{2k+1, 2k+3})|+|L_0(M_{2k+1, 2k+3})|=\sum_{i=0}^k{1\over i+1}{2i\choose i} {2k-2i\choose k-i}
$$
as desired, which completes the proof. \qed

Recall that $J_i(P_{t,t+1})$ is the set of order ideals  in $J(P_{t,t+1})$ such that the $i$ is the smallest missing  element of rank $0$.
  For $t\geq 2$ and $2\leq i\leq t$, denote by  $J^*_i(P_{t,t+1})$   the set  of order ideals  in $J_i(P_{t,t+1})$ in which exactly one element less than $i$ is marked by a star.
Let $J^*(P_{t,t+1})=\bigcup_{i=2}^{t}J^*_i(P_{t,t+1})$.

\begin{lemma}\label{lem1.4}
For $t\geq 1$, we have $|J^*(P_{t+1,t+2})|= \sum_{i=1}^{t}ic_ic_{t-i}$.
\end{lemma}
\pf From the construction of the map $\phi$, it is easily seen the the map $\phi$ sends an order ideal in  $J^*(P_{t+1,t+2})$ to a Dyck path from $(0,0)$ to $(2t+2, 0)$ in which   exactly  one up step which is to  the left of the first return point and to the right of the first up step is marked by a star. Hence, we have
$$
|J^*(P_{t+1,t+2})|=\sum_{i=1}^{t}ic_ic_{t-i},
$$ as desired. This completes the proof. \qed

For $k\geq 1$, denote by $S(P_{2k+1, 2k+3})$ the set of  nice order ideals in $J(P_{2k+1, 2k+3})$ which contain the element  $2k+2$. For $1\leq m \leq k$ and $l\leq k-m$, denote by $S_{m,\ell}(P_{2k+1, 2k+3})$ the set of all order ideals    $I\in S(P_{2k+1, 2k+3})$ such that  $2m+1$ is the smallest missing odd element from $I$  and   $2m+2\ell$ is the largest missing even element of rank $0$.

Before we deal with the enumeration of the order ideals in $S(P_{2k+1, 2k+3})$, we need the following lemma.

\begin{lemma}\label{lem1.5}
For $k\geq 1$, we have $\sum_{i=0}^{k-1}c_i  {2k-2i-1\choose k-i} ={2k\choose k}-c_k$
\end{lemma}
\pf  Notice that  the number of free Dyck paths from $(0,0)$ to $(2k,0)$ whose steps are weakly above the $x$-axis is given by $c_k$. Hence the number of free Dyck paths  from $(0,0)$ to $(2k,0)$ which have at least one step below $x$-axis is given by ${2k\choose k}-c_k$. On the other hand,  let   $P=p_1p_2\ldots p_{2k}\in \mathcal{DF}_k$ where each $p_i\in \{U,D\}$. If $P$ has at least one step below the $x$-axis, then  suppose that $p_{2i+1}$ is the leftmost down step which is weakly below the $x$-axis. Then the section $p_1p_2\ldots p_{2i}$ is a Dyck path from $(0,0)$ to $(2i,0)$, and the remaining section of $P$ is free Dyck path from $(2i,0)$ to $(2k,0)$ which stars with a down step. Hence, the number of free Dyck paths from $(0,0)$ to $(2k,0)$ which have at least one down step weakly below the $x$-axis is given by $\sum_{i=0}^{k-1}c_i{2k-2i-1\choose k-i}$.  Hence, we have
$$
\sum_{i=0}^{k-1}c_i  {2k-2i-1\choose k-i}={2k\choose k}-c_k
$$
as desired. This completes the proof. \qed

Now we are in the position to get the enumeration of nice order ideals of $J(P_{2k+1, 2k+3})$ which  contain the element $2k+2$.
Let $X_{k,\ell}$ denot the set of ordered pairs  $(I', I'')$ where  $I'\in \widetilde{L}(M_{2\ell+1, 2\ell+3})$ and $ I''\in J^*(P_{k-\ell+1, k-\ell+2})  $. Let $X_k=\bigcup_{\ell=0}^{k-1} X_{k,\ell}$.

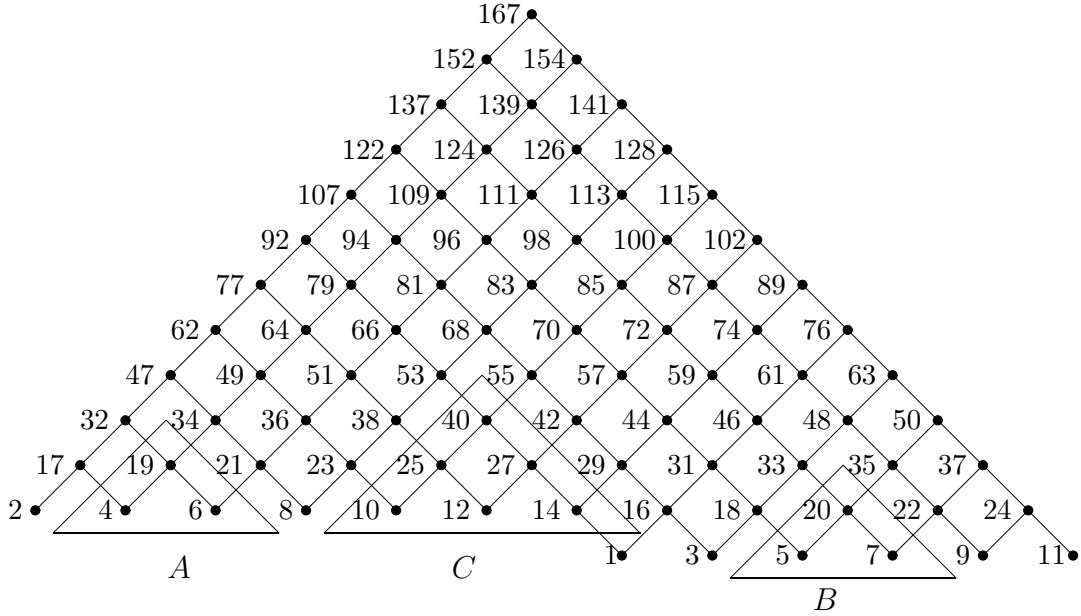
\begin{figure}[h,t]
\begin{center}
\begin{picture}(120,80)
\setlength{\unitlength}{6mm}

 \put(0,2){\circle*{0.2}}
 \put(0,2){\line(1,1){1}}

 \put(1,3){\circle*{0.2}}
 \put(1,3){\line(1,-1){1}}
 \put(2,2){\circle*{0.2}}
 \put(2,2){\line(1,1){1}}
 \put(3,3){\circle*{0.2}}
 \put(3,3){\line(1,-1){1}}

 \put(4,2){\circle*{0.2}}
 \put(4,2){\line(1,1){1}}
 \put(5,3){\circle*{0.2}}
 \put(5,3){\line(1,-1){1}}

 \put(6,2){\circle*{0.2}}
 \put(6,2){\line(1,1){1}}\put(7,3){\circle*{0.2}}
 \put(7,3){\line(1,-1){1}}\put(8,2){\circle*{0.2}}
 \put(8,2){\line(1, 1){1}}\put(9,3){\circle*{0.2}}
 \put(9,3){\line(1,-1){1}}\put(10,2){\circle*{0.2}}
 \put(10,2){\line(1,1){1}}\put(11,3){\circle*{0.2}}
 \put(11,3){\line(1,-1){1}}\put(12,2){\circle*{0.2}}
 \put(12,2){\line(1,1){1}}\put(13,3){\circle*{0.2}}
 \put(13,3){\line(1,-1){1}}\put(14,2){\circle*{0.2}}
 \put(14,2){\line(1,1){1}}\put(15,3){\circle*{0.2}}
 \put(15,3){\line(1,-1){1}}\put(16,2){\circle*{0.2}}
 \put(16,2){\line(1, 1){1}}\put(17,3){\circle*{0.2}}
 \put(17,3){\line(1, -1){1}}\put(18,2){\circle*{0.2}}
 \put(18,2){\line(1, 1){1}}\put(19,3){\circle*{0.2}}
 \put(19,3){\line(1, -1){1}}\put(20,2){\circle*{0.2}}
 \put(20,2){\line(1, 1){1}}\put(21,3){\circle*{0.2}}
 \put(21,3){\line(1, -1){1}}\put(22,2){\circle*{0.2}}
   \put(2,4){\circle*{0.2}}\put(2,4){\line(-1, -1){1}}
 \put(2,4){\line(1, -1){1}}

 \put(4,4){\circle*{0.2}}\put(4,4){\line(-1, -1){1}}
 \put(4,4){\line(1, -1){1}}

  \put(6,4){\circle*{0.2}}\put(6,4){\line(-1, -1){1}}
 \put(6,4){\line(1, -1){1}}
 \put(8,4){\circle*{0.2}}\put(8,4){\line(-1, -1){1}}
 \put(8,4){\line(1, -1){1}}
 \put(10,4){\circle*{0.2}}\put(10,4){\line(-1, -1){1}}
 \put(10,4){\line(1, -1){1}}
 \put(12,4){\circle*{0.2}}\put(12,4){\line(-1, -1){1}}
 \put(12,4){\line(1, -1){1}}
 \put(14,4){\circle*{0.2}}\put(14,4){\line(-1, -1){1}}
 \put(14,4){\line(1, -1){1}}
 \put(16,4){\circle*{0.2}}\put(16,4){\line(-1, -1){1}}
 \put(16,4){\line(1, -1){1}}
 \put(18,4){\circle*{0.2}}\put(18,4){\line(-1, -1){1}}
 \put(18,4){\line(1, -1){1}}
 \put(20,4){\circle*{0.2}}\put(20,4){\line(-1, -1){1}}
 \put(20,4){\line(1, -1){1}}

 \put(3,5){\circle*{0.2}}\put(3,5){\line(-1, -1){1}}
 \put(3,5){\line(1, -1){1}}

 \put(5,5){\circle*{0.2}}\put(5,5){\line(-1, -1){1}}
 \put(5,5){\line(1, -1){1}}

  \put(7,5){\circle*{0.2}}\put(7,5){\line(-1, -1){1}}
 \put(7,5){\line(1, -1){1}}
 \put(9,5){\circle*{0.2}}\put(9,5){\line(-1, -1){1}}
 \put(9,5){\line(1, -1){1}}
 \put(11,5){\circle*{0.2}}\put(11,5){\line(-1, -1){1}}
 \put(11,5){\line(1, -1){1}}
 \put(13,5){\circle*{0.2}}\put(13,5){\line(-1, -1){1}}
 \put(13,5){\line(1, -1){1}}
 \put(15,5){\circle*{0.2}}\put(15,5){\line(-1, -1){1}}
 \put(15,5){\line(1, -1){1}}
 \put(17,5){\circle*{0.2}}\put(17,5){\line(-1, -1){1}}
 \put(17,5){\line(1, -1){1}}
 \put(19,5){\circle*{0.2}}\put(19,5){\line(-1, -1){1}}
 \put(19,5){\line(1, -1){1}}

 \put(4,6){\circle*{0.2}}\put(4,6){\line(-1, -1){1}}
 \put(4,6){\line(1, -1){1}}

 \put(6,6){\circle*{0.2}}\put(6,6){\line(-1, -1){1}}
 \put(6,6){\line(1, -1){1}}

  \put(8,6){\circle*{0.2}}\put(8,6){\line(-1, -1){1}}
 \put(8,6){\line(1, -1){1}}
 \put(10,6){\circle*{0.2}}\put(10,6){\line(-1, -1){1}}
 \put(10,6){\line(1, -1){1}}
 \put(12,6){\circle*{0.2}}\put(12,6){\line(-1, -1){1}}
 \put(12,6){\line(1, -1){1}}
 \put(14,6){\circle*{0.2}}\put(14,6){\line(-1, -1){1}}
 \put(14,6){\line(1, -1){1}}
 \put(16,6){\circle*{0.2}}\put(16,6){\line(-1, -1){1}}
 \put(16,6){\line(1, -1){1}}
 \put(18,6){\circle*{0.2}}\put(18,6){\line(-1, -1){1}}
 \put(18,6){\line(1, -1){1}}

  \put(5,7){\circle*{0.2}}\put(5,7){\line(-1, -1){1}}
 \put(5,7){\line(1, -1){1}}

 \put(7,7){\circle*{0.2}}\put(7,7){\line(-1, -1){1}}
 \put(7,7){\line(1, -1){1}}

  \put(9,7){\circle*{0.2}}\put(9,7){\line(-1, -1){1}}
 \put(9,7){\line(1, -1){1}}
 \put(11,7){\circle*{0.2}}\put(11,7){\line(-1, -1){1}}
 \put(11,7){\line(1, -1){1}}
 \put(13,7){\circle*{0.2}}\put(13,7){\line(-1, -1){1}}
 \put(13,7){\line(1, -1){1}}
 \put(15,7){\circle*{0.2}}\put(15,7){\line(-1, -1){1}}
 \put(15,7){\line(1, -1){1}}
 \put(17,7){\circle*{0.2}}\put(17,7){\line(-1, -1){1}}
 \put(17,7){\line(1, -1){1}}

 \put(6,8){\circle*{0.2}}\put(6,8){\line(-1, -1){1}}
 \put(6,8){\line(1, -1){1}}

 \put(8,8){\circle*{0.2}}\put(8,8){\line(-1, -1){1}}
 \put(8,8){\line(1, -1){1}}

  \put(10,8){\circle*{0.2}}\put(10,8){\line(-1, -1){1}}
 \put(10,8){\line(1, -1){1}}
 \put(12,8){\circle*{0.2}}\put(12,8){\line(-1, -1){1}}
 \put(12,8){\line(1, -1){1}}
 \put(14,8){\circle*{0.2}}\put(14,8){\line(-1, -1){1}}
 \put(14,8){\line(1, -1){1}}
 \put(16,8){\circle*{0.2}}\put(16,8){\line(-1, -1){1}}
 \put(16,8){\line(1, -1){1}}

 \put(7,9){\circle*{0.2}}\put(7,9){\line(-1, -1){1}}
 \put(7,9){\line(1, -1){1}}

 \put(9,9){\circle*{0.2}}\put(9,9){\line(-1, -1){1}}
 \put(9,9){\line(1, -1){1}}

  \put(11,9){\circle*{0.2}}\put(11,9){\line(-1, -1){1}}
 \put(11,9){\line(1, -1){1}}
 \put(13,9){\circle*{0.2}}\put(13,9){\line(-1, -1){1}}
 \put(13,9){\line(1, -1){1}}
 \put(15,9){\circle*{0.2}}\put(15,9){\line(-1, -1){1}}
 \put(15,9){\line(1, -1){1}}

 \put(8,10){\circle*{0.2}}\put(8,10){\line(-1, -1){1}}
 \put(8,10){\line(1, -1){1}}

 \put(10,10){\circle*{0.2}}\put(10,10){\line(-1, -1){1}}
 \put(10,10){\line(1, -1){1}}

  \put(12,10){\circle*{0.2}}\put(12,10){\line(-1, -1){1}}
 \put(12,10){\line(1, -1){1}}
 \put(14,10){\circle*{0.2}}\put(14,10){\line(-1, -1){1}}
 \put(14,10){\line(1, -1){1}}

  \put(9,11){\circle*{0.2}}\put(9,11){\line(-1, -1){1}}
 \put(9,11){\line(1, -1){1}}

 \put(11,11){\circle*{0.2}}\put(11,11){\line(-1, -1){1}}
 \put(11,11){\line(1, -1){1}}

  \put(13,11){\circle*{0.2}}\put(13,11){\line(-1, -1){1}}
 \put(13,11){\line(1, -1){1}}

 \put(10,12){\circle*{0.2}}\put(10,12){\line(-1, -1){1}}
 \put(10,12){\line(1, -1){1}}

 \put(12,12){\circle*{0.2}}\put(12,12){\line(-1, -1){1}}
 \put(12,12){\line(1, -1){1}}

  \put(11,13){\circle*{0.2}}\put(11,13){\line(-1, -1){1}}
 \put(11,13){\line(1, -1){1}}

  \put(12,2){\line(1, -1){1}}\put(13,1){\circle*{0.2}}\put(13,1){\line(1, 1){1}}
  \put(15,1){\line(-1, 1){1}}\put(15,1){\circle*{0.2}}\put(15,1){\line(1, 1){1}}\put(17,1){\line(-1, 1){1}}\put(17,1){\circle*{0.2}}\put(17,1){\line(1, 1){1}}\put(19,1){\line(-1, 1){1}}\put(19,1){\circle*{0.2}}\put(19,1){\line(1, 1){1}}\put(21,1){\line(-1, 1){1}}\put(21,1){\circle*{0.2}}\put(21,1){\line(1, 1){1}} \put(23,1){\line(-1, 1){1}}\put(23,1){\circle*{0.2}}

 \put(-0.6, 1.8){\small$2$}\put(1.4, 1.8){\small$4$}\put(3.4, 1.8){\small$6$}
 \put(5.4, 1.8){\small$8$} \put(7, 1.8){\small$10$}\put(9, 1.8){\small$12$}
 \put(11, 1.8){\small$14$}\put(13, 1.8){\small$16$}\put(15, 1.8){\small$18$}
 \put(17, 1.8){\small$20$}\put(19, 1.8){\small$22$} \put(21, 1.8){\small$24$}

 \put(0, 2.8){\small$17$}\put(2, 2.8){\small$19$}\put(4, 2.8){\small$21$}
 \put(6, 2.8){\small$23$} \put(8, 2.8){\small$25$}\put(10, 2.8){\small$27$}
 \put(12, 2.8){\small$29$}\put(14, 2.8){\small$31$}\put(16, 2.8){\small$33$}
 \put(18, 2.8){\small$35$}\put(20, 2.8){\small$37$}

 \put(1, 3.8){\small$32$}\put(3, 3.8){\small$34$}\put(5, 3.8){\small$36$}
 \put(7, 3.8){\small$38$} \put(9, 3.8){\small$40$}\put(11, 3.8){\small$42$}
 \put(13, 3.8){\small$44$}\put(15, 3.8){\small$46$}\put(17, 3.8){\small$48$}
 \put(19, 3.8){\small$50$}

 \put(2, 4.8){\small$47$}\put(4, 4.8){\small$49$}\put(6, 4.8){\small$51$}
 \put(8, 4.8){\small$53$} \put(10, 4.8){\small$55$}\put(12, 4.8){\small$57$}
 \put(14, 4.8){\small$59$}\put(16, 4.8){\small$61$}\put(18, 4.8){\small$63$}

  \put(3, 5.8){\small$62$}\put(5, 5.8){\small$64$}\put(7, 5.8){\small$66$}
 \put(9, 5.8){\small$68$} \put(11, 5.8){\small$70$}\put(13, 5.8){\small$72$}
 \put(15, 5.8){\small$74$}\put(17, 5.8){\small$76$}

 \put(4, 6.8){\small$77$}\put(6, 6.8){\small$79$}\put(8, 6.8){\small$81$}
 \put(10, 6.8){\small$83$} \put(12, 6.8){\small$85$}\put(14, 6.8){\small$87$}
 \put(16, 6.8){\small$89$}

 \put(5, 7.8){\small$92$}\put(6.8, 7.8){\small$94$}\put(8.8, 7.8){\small$96$}
 \put(10.8, 7.8){\small$98$} \put(12.8, 7.8){\small$100$}\put(14.8, 7.8){\small$102$}

  \put(5.8, 8.8){\small$107$}\put(7.8, 8.8){\small$109$}\put(9.8, 8.8){\small$111$}
 \put(11.8, 8.8){\small$113$} \put(13.8, 8.8){\small$115$}

 \put(6.8, 9.8){\small$122$}\put(8.8, 9.8){\small$124$}\put(10.8, 9.8){\small$126$}
 \put(12.8, 9.8){\small$128$}

 \put(7.8, 10.8){\small$137$}\put(9.8, 10.8){\small$139$}\put(11.8, 10.8){\small$141$}

  \put(8.8, 11.8){\small$152$}\put(10.8, 11.8){\small$154$}
  \put(9.8, 12.8){\small$167$}

  \put(12.6, 0.8){\small$1$}\put(14.4, 0.8){\small$3$}\put(16.4, 0.8){\small$5$}\put(18.4, 0.8){\small$7$}\put(20.4, 0.8){\small$9$}
  \put(22.2, 0.8){\small$11$}


  \put(0.4,1.5){\line(1, 1){2.5}}\put(2.9,4){\line(1, -1){2.5}}
  \put(0.4,1.5){\line(1, 0){5}}


  \put(15.4,0.5){\line(1, 1){2.5}}\put(17.9,3){\line(1, -1){2.5}}
  \put(15.4,0.5){\line(1, 0){5}}

  \put(6.4,1.5){\line(1, 1){3.5}}\put(9.9,5){\line(1, -1){3.5}}
  \put(6.4,1.5){\line(1, 0){7}}

  \put(2.7, 0.5){  $A$  } \put(9, 0.5){  $C$  } \put(17, -0.2){  $B$  }
 \end{picture}
\end{center}
\caption{ The decomposition of an  order ideal of $S_{1,3}(P_{13,15})$.}\label{decomposition}
\end{figure}

\begin{theorem}\label{th5}
Let $k$ be a  positive integer. The cardinality of the set $S(P_{2k+1, 2k+3})$ is given by $\sum_{j=1}^k{j\over j+1}{2j\choose j} {2k-2j\choose k-j}$.
\end{theorem}
\pf Obviously, the assertion holds for $k=1$ since there is exactly one order ideal  in $S(P_{3, 5})$.  For $k\geq 2$ ,   we first describe a map $\gamma$ from the set $S(P_{2k+1, 2k+3})$  to the set $X_{k}$.
  Let $I$ be an order ideal in $S_{m,\ell}(P_{2k+1, 2k+3})$ for some  $1\leq m\leq k$ and $0\leq \ell\leq k-m$. Then  we can  decompose $I$ into three parts: one is
  $\{1, 3, \ldots, 2m-1\}$, one is isomorphic to an order ideal $I_1\in \widetilde{L}(M_{2\ell+1, 2\ell+3})$, and one is isomorphic to an order ideal $I_2\in J(P_{k-\ell+1, k-\ell+2})$.  For example, Figure \ref{decomposition} illustrates the decomposition of   an order ideal $I$ of $S_{1,3}(P_{13, 15})$. It  can be decomposed into three parts: one is $\{1\}$, one is order-isomorphic to an order ideal in $\widetilde{L}(P_{7,9})$ whose elements lie in the  triangles  $A$ and $B$, and one is order-isomorphic to an order ideal in $J(P_{4,5})$ whose elements lie in the  triangle  $C$.

  By assigning a star to the element $k-m-\ell+1$, we can obtain a marked order ideal $I'_2\in   J^*(P_{k-\ell+1, k-\ell+2})$ from $I_2$.
 Set $\gamma(I)=(I_1, I'_2)$. It is easily seen that the map $\gamma$ is well defined.  Conversely, given an ordered pair $(I', I'')\in X_{k,\ell}$ such that the element $k-\ell-m+1$ is marked by a star, we can recover the order ideal $I$ by reversing the above procedure. Hence, the map $\gamma$ is a bijection between the set $ S(P_{2k+1, 2k+3}) $ and the set $X_k$. This yields that
  $|S(P_{2k+1, 2k+3})|=|X_k|$.

According to the definition of $X_k$, we have
  
  $$
  |X_k|=\sum_{\ell=0}^{k-1}|X_{k,\ell}|
  $$
  for all $k\geq 2$.
  Notice that  the cardinality of the set $\widetilde{L}(M_{1,3})$ is given by $1$. 
  By Lemmas \ref{lem1.3} and \ref{lem1.4}, we deduce that, for $k\geq 2$, 
  $$
  \begin{array}{lll}
  |X_k|&=&\sum_{i=1}^{k}  ic_ic_{k-i}+\sum_{\ell=1}^{k-1}{2\ell-1\choose \ell}(\sum_{i=1}^{k-\ell}  ic_ic_{k-\ell-i}   )
  \\[3pt]
  &=&\sum_{i=1}^{k}  ic_ic_{k-i}+\sum_{i=1}^{k-1} ic_i(\sum_{\ell=1}^{k-i}{2\ell-1\choose \ell}c_{k-i-\ell})
  \\[3pt]
     &=& \sum_{i=1}^{k}  ic_ic_{k-i}+\sum_{i=1}^{k-1} ic_i(\sum_{m=0}^{k-i-1}c_m{2k-2i-2m-1\choose k-i-m})  
    \\[3pt]
       &=&\sum_{i=1}^{k}  ic_ic_{k-i}+\sum_{i=1}^{k-1}ic_i ({2k-2i\choose k-i}-c_{k-i} ) \,\,\,\,\, \mbox{(by Lemma \ref{lem1.5})}\\[5pt]
  &=& \sum_{i=1}^k{i\over i+1}{2i\choose i}  {2k-2i\choose k-i}
  \end{array}
  $$ as desired. This completes the proof.\qed

We are now ready to complete the proof of Conjecture \ref{con1}.

\noindent
{\it Proof of Conjecture \ref{con1}.}
Combining Theorems \ref{th1.4}  and \ref{th5},  we deduce  that the number of $(2k+1, 2k+3)$-core partitions into distinct parts is given by
$$|L(M_{2k+1, 2k+3})|+|S(P_{2k+1, 2k+3})|=\sum_{i=0}^{k}{2i\choose i}{2k-2i
\choose k-i}.
$$
   Notice that
$${1\over 1-4x}=({1\over \sqrt{1-4x}})^2 =(\sum_{n\geq 0}{2n\choose n}x^n)^2.$$  By comparing the coefficient of $x^k$ in the above formula,  we have $\sum_{i=0}^{k}{2i\choose i}{2k-2i
\choose k-i}=2^{2k}$. This yields that
\begin{equation}\label{eq1}
|L(M_{2k+1, 2k+3})|+|S(P_{2k+1, 2k+3})|=2^{2k}.
 \end{equation}
 Substituting  $k={s-1\over 2}$ in (\ref{eq1}), we are led to a proof of   Conjecture \ref{con1}.
 \qed

\section{Proof of Conjecture \ref{con2}}
In this section, we shall construct a partition $\kappa_k$ for any integer $k\geq 0$. It turns out that $\kappa_k$ is the unique $(2k+1, 2k+3)$-core partitions into distinct parts which has the largest size.   This leads to a proof of Conjecture \ref{con2}.

The following  four lemmas give a characterization of the  largest $(2k+1, 2k+3)$-core partitions  with distinct parts whose  $\beta$-set is contained in $L(M_{2k+1, 2k+3})$.
Denote by $B_s$ (resp. $B'_s$ )  the set of elements of rank $s$ in $T_k$ (resp. $T'_k$).
By Lemma \ref{mainlem2}, we have
$B_s=\{2i-1+(s-1)(2k+3)\mid  1\leq i\leq k+1-s\}$ and $B'_s=\{2i+(s-1)(2k+3)\mid  1\leq i\leq k+1-s\}$ for all $1\leq s\leq k$.

\begin{lemma}\label{2.1}
Let $\lambda$ be a $(2k+1, 2k+3)$-core partitions into distinct parts which has the largest size.     If $\beta(\lambda)\in L(M_{2k+1,2k+3})$ and $\beta(\lambda)$ contains an element $x$ in  $B_s$(resp. $B'_s$), then $\beta(\lambda)$      contains       all the element $y\in B_s$  ( resp. $y\in B'_s$)  when $y> x$.
\end{lemma}

\pf   It suffices to show that  $\beta(\lambda)$   also  contains    the element $x+2$ if    $\beta(\lambda)$ contains an element $x$ in  $B_s$(resp. $B'_s$) and $x+2\in B_s $ (resp. $x+2\in B'_s $).    If not, we choose $s$ to be the smallest integer for such  $B_s$ (resp. $B'_s$) and let $x$ be the smallest such number once $s$ is determined. If $s\geq 1$, then replace  $y$ by $y+2$ if $y\succeq p$ for some $p\leq x$ and $p\in B_s$ (resp. $p\in B'_s$).   Otherwise,  replace $y$ by $y+1$ if $y\succeq p$ for some $p\leq x$ and $p\in B_0$ (resp. $p\in B'_0$). By this process, we obtain a new nice order ideal  $\beta'$ of $M_{2k+1, 2k+3}$. It is easily seen that $\beta'$ has the same number of parts as $\beta(\lambda)$ and  a larger sum of the elements. By relation (\ref{eq}), the size of the $(2k+1, 2k+3)$-core partition into distinct parts corresponding to $\beta'$ is larger than $\lambda$, which contradicts the assumption that $\lambda$ is a $(2k+1, 2k+3)$-core partitions into distinct parts which has the largest size. This completes the proof. \qed

 Denote by $U_t$ the unique order ideal which contains all the elements of $P_{t,t+1}$.
For $1\leq i\leq 2k$, let $\beta_{i,0}$   denote the unique nice order ideal in  $M_{2k+1, 2k+3}$  which is isomorphic to $U_{k-\lfloor{i+1\over 2}\rfloor+1}$ and contains all the elements $ i+2,  i+4, \ldots, i+ 2(k-\lfloor{i+1\over 2}\rfloor)   $.  For $1\leq i\leq 2k$ and $1\leq j\leq k- \lfloor{i+1\over 2}\rfloor+1$, let $\beta_{i,j}$ be the union of $\beta_{i,0}$ and the chain consisting of $i, i+2k+3, \ldots, i+(j-1)(2k+3)$. See Figure  \ref{beta{2,3}} for an illustration of the order ideal $\beta_{4,3}$.  For $1\leq j\leq k$, denote by $\gamma_{k,j}$ be the union of $\beta_{1,k}$ and the chain consisting of $2k+2, 2k+2+(2k+3), \ldots, 2k+2+(j-1)(2k+3)$.

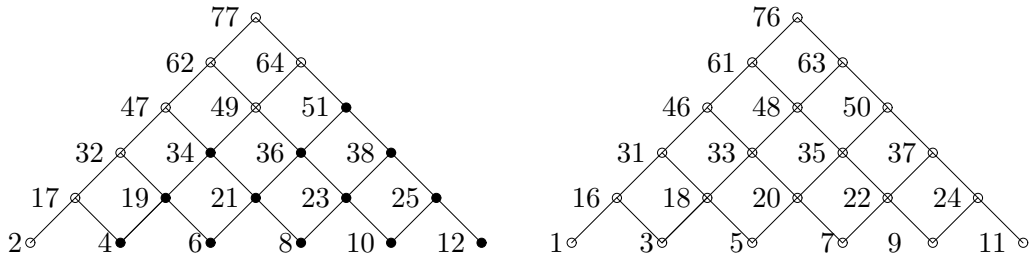
\begin{figure}[h,t]
\begin{center}
\begin{picture}(120,40)
\setlength{\unitlength}{6mm}

 \put(0,0){\circle {0.2}}
 \put(0,0){\line(1,1){1}}\put(1,1){\circle {0.2}}
 \put(1,1){\line(1,-1){1}}\put(2,0){\circle* {0.2}}
 \put(2,0){\line(1,1){1}} \put(2,0){\line(1,1){1}}\put(3,1){\circle* {0.2}}
 \put(3,1){\line(1,-1){1}}\put(4,0){\circle* {0.2}}
 \put(4,0){\line(1,1){1}}\put(5,1){\circle*{0.2}}
 \put(5,1){\line(1,-1){1}}\put(6,0){\circle*{0.2}}
 \put(6,0){\line(1,1){1}}\put(7,1){\circle*{0.2}}
 \put(7,1){\line(1,-1){1}}\put(8,0){\circle*{0.2}}
 \put(8,0){\line(1,1){1}}\put(9,1){\circle*{0.2}}
 \put(9,1){\line(1,-1){1}}\put(10,0){\circle*{0.2}}

 \put(2,2){\circle{0.2}}\put(2,2){\line(1,-1){1}}
 \put(2,2){\line(-1,-1){1}}\put(4,2){\circle*{0.2}}\put(4,2){\line(1,-1){1}}
 \put(4,2){\line(-1,-1){1}}\put(6,2){\circle*{0.2}}\put(6,2){\line(1,-1){1}}
 \put(6,2){\line(-1,-1){1}}\put(8,2){\circle*{0.2}}\put(8,2){\line(1,-1){1}}
 \put(8,2){\line(-1,-1){1}}

 \put(3,3){\circle{0.2}}\put(3,3){\line(1,-1){1}}
 \put(3,3){\line(-1,-1){1}}\put(5,3){\circle{0.2}}\put(5,3){\line(1,-1){1}}
 \put(5,3){\line(-1,-1){1}}\put(7,3){\circle*{0.2}}\put(7,3){\line(1,-1){1}}
 \put(7,3){\line(-1,-1){1}}

 \put(4,4){\circle{0.2}}\put(4,4){\line(1,-1){1}}
 \put(4,4){\line(-1,-1){1}}\put(6,4){\circle{0.2}}\put(6,4){\line(1,-1){1}}
 \put(6,4){\line(-1,-1){1}}

 \put(5,5){\circle{0.2}}\put(5,5){\line(1,-1){1}}
 \put(5,5){\line(-1,-1){1}}
 \put(-0.5, -0.2){\small$2$}\put(1.5, -0.2){\small$4$}\put(3.5, -0.2){\small$6$}\put(5.5, -0.2){\small$8$}\put(7, -0.2){\small$10$}
 \put(9, -0.2){\small$12$}

 \put(0, 0.8){\small$17$}\put(2, 0.8){\small$19$}\put(4, 0.8){\small$21$}\put(6, 0.8){\small$23$}\put(8, 0.8){\small$25$}

  \put(1, 1.8){\small$32$}\put(3, 1.8){\small$34$}\put(5, 1.8){\small$36$}\put(7, 1.8){\small$38$}

  \put(2, 2.8){\small$47$}\put(4, 2.8){\small$49$}\put(6, 2.8){\small$51$}

  \put(3, 3.8){\small$62$}\put(5, 3.8){\small$64$}
  \put(4, 4.8){\small$77$}


  \put(12,0){\circle{0.2}}
 \put(12,0){\line(1,1){1}}\put(13,1){\circle{0.2}}
 \put(13,1){\line(1,-1){1}}\put(14,0){\circle{0.2}}
 \put(14,0){\line(1,1){1}} \put(14,0){\line(1,1){1}}\put(15,1){\circle{0.2}}
 \put(15,1){\line(1,-1){1}}\put(16,0){\circle{0.2}}
 \put(16,0){\line(1,1){1}}\put(17,1){\circle{0.2}}
 \put(17,1){\line(1,-1){1}}\put(18,0){\circle{0.2}}
 \put(18,0){\line(1,1){1}}\put(19,1){\circle{0.2}}
 \put(19,1){\line(1,-1){1}}\put(20,0){\circle{0.2}}
 \put(20,0){\line(1,1){1}}\put(21,1){\circle{0.2}}
 \put(21,1){\line(1,-1){1}}\put(22,0){\circle{0.2}}

 \put(14,2){\circle{0.2}}\put(14,2){\line(1,-1){1}}
 \put(14,2){\line(-1,-1){1}}\put(16,2){\circle{0.2}}\put(16,2){\line(1,-1){1}}
 \put(16,2){\line(-1,-1){1}}\put(18,2){\circle{0.2}}\put(18,2){\line(1,-1){1}}
 \put(18,2){\line(-1,-1){1}}\put(20,2){\circle{0.2}}\put(20,2){\line(1,-1){1}}
 \put(20,2){\line(-1,-1){1}}

 \put(15,3){\circle{0.2}}\put(15,3){\line(1,-1){1}}
 \put(15,3){\line(-1,-1){1}}\put(17,3){\circle{0.2}}\put(17,3){\line(1,-1){1}}
 \put(17,3){\line(-1,-1){1}}\put(19,3){\circle{0.2}}\put(19,3){\line(1,-1){1}}
 \put(19,3){\line(-1,-1){1}}

 \put(16,4){\circle{0.2}}\put(16,4){\line(1,-1){1}}
 \put(16,4){\line(-1,-1){1}}\put(18,4){\circle{0.2}}\put(18,4){\line(1,-1){1}}
 \put(18,4){\line(-1,-1){1}}

 \put(17,5){\circle{0.2}}\put(17,5){\line(1,-1){1}}
 \put(17,5){\line(-1,-1){1}}
 \put(11.5, -0.2){\small$1$}\put(13.5, -0.2){\small$3$}\put(15.5, -0.2){\small$5$}\put(17.5, -0.2){\small$7$}\put(19, -0.2){\small$9$}
 \put(21, -0.2){\small$11$}

 \put(12, 0.8){\small$16$}\put(14, 0.8){\small$18$}\put(16, 0.8){\small$20$}\put(18, 0.8){\small$22$}\put(20, 0.8){\small$24$}

  \put(13, 1.8){\small$31$}\put(15, 1.8){\small$33$}\put(17, 1.8){\small$35$}\put(19, 1.8){\small$37$}

  \put(14, 2.8){\small$46$}\put(16, 2.8){\small$48$}\put(18, 2.8){\small$50$}

  \put(15, 3.8){\small$61$}\put(17, 3.8){\small$63$}
  \put(16, 4.8){\small$76$}
\end{picture}
\end{center}
\caption{ The order ideal  $\beta_{4,3}$.}\label{beta{2,3}}
\end{figure}

  It is easy to check that both $\beta_{i,j}$ and $\gamma_{k,j}$ are nice order ideals of $P_{2k+1, 2k+3}$.  Let $\lambda_{i,j}$ (resp. $\mu_{ k,j}$) be the unique partition such that $\beta(\lambda_{i,j})=\beta_{i,j}$ (resp. $\beta(\mu_{k,j})=\gamma_{k,j}$ ).   By Theorem \ref{th1.3},  both $\lambda_{ i,j}$ and $\mu_{k,j}$ are $(2k+1, 2k+3)$-core partitions with distinct parts.

\begin{lemma}\label{lem2.2}
Fix $k\geq 1$. Let $\lambda$ be a $(2k+1, 2k+3)$-core partition into distinct parts such that $\beta(\lambda)\in L(M_{2k+1, 2k+3})$.  If $\lambda$ if of maximum size, then there exist some integers $1\leq i\leq 2k$ and $1\leq j\leq k- \lfloor{i+1\over 2}\rfloor+1$  such that  $\beta(\lambda)=\beta_{i,j}$.
\end{lemma}

\pf Let $i$ be the minimal integer of rank $0$   in $\beta(\lambda)$ and $j$  be the maximal integer such that $i+(j-1)(2k+3)$ is contained in $\beta(\lambda)$. We claim that $i+2+m(2k+3)$ is contained in $\beta(\lambda)$ for all $0\leq m\leq k- \lfloor{i+1\over 2}\rfloor-1$. Assume to the contrary that the claim is not valid, that is,  there exists some integer $\ell$ such that $\beta(\lambda)$ does not contain the element $i+2+\ell(2k+3)$.
  We choose $s$ to be the smallest such integer.  By Lemma \ref{2.1}, we have $s> j-1$. By replacing the element $i+(j-1)(2k+3)$ by the element $i+2+s(2k+3)$, we get a new nice order ideal $\beta'$  of $J(P_{2k+1, 2k+3})$. It is easy to check that $\beta'$ has the same cardinality as $\beta(\lambda)$ and has a larger sum of the elements. By relation (\ref{eq}), the size of the partition corresponding to $\beta'$ is larger than that of $\lambda$. This yields a contradiction with the fact that $\lambda$ has the largest size, which completes the proof of the claim. Combining the claim and Lemma \ref{2.1}, we have $\beta(\lambda)=\beta_{i,j}$. This completes the proof. \qed

\begin{lemma}\label{lem2.4}
For $1\leq i\leq 2k$, we have $|\lambda_{i,j}|\leq |\lambda_{i, k- \lfloor{i+1\over 2}\rfloor+1}|$ for $1\leq j\leq k- \lfloor{i+1\over 2}\rfloor+1$, with the equality holding if and only if $j=k- \lfloor{i+1\over 2}\rfloor+1$. Moreover, we have $|\lambda_{i,k- \lfloor{i+1\over 2}\rfloor+1}|>|\lambda_{i+2, k- \lfloor{i+1\over 2}\rfloor}|$ for all $k\geq 3$ and  $1\leq i\leq 2k-2$.
\end{lemma}

\pf By relation (\ref{eq}), the size of $\lambda_{i,j}$ is given by
\begin{equation}\label{eq2.1}
\begin{array}{lll}
|\lambda_{i,j}|&=& \sum_{h\in \beta_{i,j}}h-{|\beta_{i,j}|\choose 2}\\[5pt]
&=& \sum_{h\in \beta_{i,0}}h+\sum_{p=0}^{j-1}(i+p(2k+3))-{|\beta_{i,0}|+j\choose 2}\\[5pt]
&=& \sum_{h\in \beta_{i,0}}h-{|\beta_{i,0}|\choose 2}+ \sum_{p=0}^{j-1}(i+p(2k+3))-{|\beta_{i,0}|+j\choose 2}+{|\beta_{i,0}|\choose 2}\\[5pt]
&=& |\lambda_{i,0}| +ij+{j\choose 2}(2k+2)-j|\beta_{i,0}|
\end{array}
\end{equation}
By the definition of $\beta_{i,0}$, we obtain that
\begin{equation}\label{eq2.2}
|\beta_{i,0}|=
{k-{\lfloor {i+1\over 2} \rfloor}+1\choose 2}.
 \end{equation}
Hence, we have
\begin{equation}\label{eq2.3}
|\lambda_{i,j}|=|\lambda_{i,0}| +ij+{j\choose 2}(2k+2)-j{k-{\lfloor {i+1\over 2} \rfloor}+1\choose 2}.
\end{equation}
In particular,   we have
\begin{equation}\label{eq2.4}
|\lambda_{i, k-{\lfloor {i+1\over 2} \rfloor}+1}|-|\lambda_{i,0}|=\left\{\begin{array}{ll}
(k+m+1){k-m+1\choose 2}+2m(k-m+1)&\,\,\, \mbox{if}\,\, $i=2m$ ,\\[5pt]
(k+m+1){k-m+1\choose 2}+(2m-1)(k-m+1)&\,\,\, \mbox{if}\,\, $i=2m-1$ ,\\[5pt]
\end{array}\right.
\end{equation}
and

 \begin{equation}\label{eq2.5}
|\lambda_{i, 1}|-|\lambda_{i,0}|=\left\{\begin{array}{ll}
 -{k-m+1\choose 2}+2m &\,\,\, \mbox{if}\,\, $i=2m$ ,\\[5pt]
 -{k-m+1\choose 2}+2m-1 &\,\,\, \mbox{if}\,\, $i=2m-1$ .\\[5pt]
\end{array}\right.
\end{equation}
This implies that
 if $i<2k-1$, then
\begin{equation}\label{eq2.6}
|\lambda_{i, k-{\lfloor {i+1\over 2} \rfloor}+1}|> |\lambda_{i,1}|.
\end{equation}
Moreover, when $i=2k $ or $2k-1$, we have
\begin{equation}\label{eq2.7}
|\lambda_{i, k-{\lfloor {i+1\over 2} \rfloor}+1}|= |\lambda_{i,1}|.
\end{equation}

For fixed $k$ and $i$, we see that $\lambda_{i,j}$ is a quadratic function of $j$ with a positive leading coefficient. Hence the maximum value of $\lambda_{i,j}$ is obtained at $j=1$ or $j= k-{\lfloor {i+1\over 2} \rfloor}+1 $ when $j$ ranges over $[1,  k-{\lfloor {i+1\over 2} \rfloor}+1 ]$. In view of (\ref{eq2.6}) and (\ref{eq2.7}), we have   $|\lambda_{i,j}|\leq |\lambda_{i, k- \lfloor{i+1\over 2}\rfloor+1}|$ for $1\leq j\leq k- \lfloor{i+1\over 2}\rfloor+1$, with the equality holding if and only if $j=k- \lfloor{i+1\over 2}\rfloor+1$.  This completes the proof of  the first part of the lemma.

From the definition of $\beta_{i+1,0}$, we have $\beta_{i, k-{\lfloor {i+1\over 2} \rfloor}+1}=\beta_{i-2, 0}$ for all $3\leq i\leq 2k$.  In view of
(\ref{eq2.4}), we have $|\lambda_{i, k-{\lfloor {i+1\over 2} \rfloor}+1}|>|\lambda_{i,0}|=|\lambda_{i+2, k-{\lfloor {i+1 \over 2} \rfloor}}|$ for all $k\geq 3$ and $1\leq i\leq 2k-2$. This completes the proof.
 \qed

\begin{lemma}\label{lem2.5}
Fix $k\geq 1$. Let $\lambda$ be a $(2k+1, 2k+3)$-core partition into distinct parts such that $\beta(\lambda)\in L(M_{2k+1, 2k+3})$. If $\lambda$ is of maximum size, then   $\beta(\lambda)=\beta_{2,k}$.
\end{lemma}
\pf From lemmas \ref{lem2.2} and \ref{lem2.4}, we have $\beta(\lambda)=\beta_{2,k}$ or $\beta_{1,k}$. According to the definition of $\beta_{1,k}$ and $\beta_{2,k}$, it is easy to verify that
  $$
  |\beta_{1,k}|=|\beta_{2,k}|
  $$
  and
   $$
   \sum_{h\in\beta_{2,k } }h=\sum_{h\in\beta_{1,k} }h+{k(k+1)\over  2}.
  $$
  By relation (\ref{eq}), we have
   \begin{equation}\label{eq2.10}
   |\lambda_{2,k}|-|\lambda_{1,k}| = {k(k+1)\over  2}>0.
 \end{equation}
 Hence we have $\beta(\lambda)=\beta_{2,k}$ as desired, which completes the proof. \qed

 Now we proceed to deal with the   largest $(2k+1, 2k+3)$-core partitions  with distinct parts  and  its $\beta$-set is contained in $S(P_{2k+1, 2k+3})$. Given positive integers $a\leq b$, we denote $\{a, a+1, \ldots, b\}$ by $[a,b]$.

Recall that $Q_k=\{2i-1\mid 1\leq i\leq k\}$ and $A_s=\{2i+(s-1)(2k+3)\mid   1\leq i\leq 2k+1-s\} $  for $1\leq s\leq 2k$.

 \begin{lemma}\label{22.1}
Let $\lambda$ be a $(2k+1, 2k+3)$-core partitions into distinct parts which has the largest size.     If $\beta(\lambda)\in S(P_{2k+1,2k+3})$, then $\beta(\lambda)$   contains        all the   elements  of  $Q_k$.
\end{lemma}
\pf     If not,   suppose that  $2m+1$ is  the smallest missing odd element  of rank $0$ from $\beta(\lambda)$ for some $1\leq m\leq k-1 $. Let  $ \ell$  be the largest integer such that the element $2m+2\ell$ is missing from $\beta(\lambda)$ once $m$ is determined. By the definition of $S_{m, \ell}(P_{2k+1, 2k+3})$, we have $\beta(\lambda)\in S_{m, \ell}(P_{2k+1, 2k+3})$. As in the proof of Theorem \ref{th5},     $\beta(\lambda)$ can be decomposed into three parts $I_1, I_2, I_3$, where
\begin{itemize}
\item  $I_1=\{1,3, \ldots, 2m-1\}$;
\item $I_2$ is isomorphic to an order ideal of $\widetilde{L}(M_{2\ell+1, 2\ell+3})$;
    \item $I_3$ is isomorphic to an order ideal of $J(P_{k-\ell+1, k-\ell+2})$.
\end{itemize}
To be more precise, $I_2$ consists of the elements
  $y\in I$  such that $y\succeq p$ for some
$p\in [2m+1, 2m+2\ell]$, whereas  $I_1\cup I_3$ consists of the elements $y\in I$  such that $y\succeq p$ for some
$p\in [1, 2m]\cup [2m+2\ell+1, 2k+2]$.

We proceed to construct a new order ideal $\beta'  $ of  $P_{2k+1, 2k+3}$ by the following procedure.
\begin{itemize}
\item Firstly, replace $I_1$ by $\{1,3, \ldots, 2m+1\}$;
\item Secondly, replace each element $y$ of $I_2$ by $y+1$;
\item Finally, replace each element $y$ of $I_3$ by $y+2$.
\end{itemize}
 Suppose that the cardinality of $I_3$ is given by $x$.
 It is easily seen that $$|\beta'|=|\beta(\lambda)|+1$$
  and $$\sum_{h\in \beta'}h=\sum_{h\in \beta(\lambda)}h+2x+2m+1+|\beta(\lambda)|-(x+m)=\sum_{h\in \beta(\lambda)}h+x+m+1+|\beta(\lambda)|.$$  By relation (\ref{eq}), the size of the $(2k+1, 2k+3)$-core partition into distinct parts corresponding to $\beta'$ is given by
  $$
  \begin{array}{lll}
  \sum_{h\in \beta'}h-{|\beta'| \choose 2}&=& \sum_{h\in \beta(\lambda)}h + |\beta(\lambda)|+x+m+1-{|\beta(\lambda)|+1\choose 2}\\
  &=& \sum_{h\in \beta(\lambda)}h-{|\beta(\lambda)| \choose 2} + |\beta(\lambda)|+1+x+m-{|\beta(\lambda)|+1\choose 2}+{|\beta(\lambda)| \choose 2}\\
  &=& |\lambda|+x+m+1.
  \end{array}
  $$
   This  yields that the partition corresponding to $\beta'$ has a larger size than  $\lambda$,
   which contradicts the assumption that $\lambda$ is a $(2k+1, 2k+3)$-core partitions into distinct parts which has the largest size. This completes the proof. \qed

 \begin{lemma}\label{22.2}
Fix $ s\geq 1$. Let $\lambda$ be a $(2k+1, 2k+3)$-core partitions into distinct parts which has the largest size.     If $\beta(\lambda)\in S(P_{2k+1,2k+3})$ and $\beta(\lambda)$ contains an element $x\in A_s$, then $\beta(\lambda)$    contains all  the elements $y\in A_s $  when $y> x$.
\end{lemma}
\pf It suffices to show that $\beta(\lambda)$   also  contains  the element $x+2$ if $\beta(\lambda)$ contains an element $x\in A_s$.   If not, we choose $s$ to be the smallest such  integer   and let $x$ be the smallest such number once $s$ is determined.
We  can obtain a new nice order ideal  $\beta'$ of $P_{2k+1, 2k+3}$  by  replacing each element   $y$ by $y+2$ if $y\succeq p$  such that  $p\leq x$ and $p$ is of rank $s$.  It is easily seen that $\beta'$ has the same number of parts as $\beta(\lambda)$ and  a larger sum of the elements. By relation (\ref{eq}), the size of the $(2k+1, 2k+3)$-core partition into distinct parts corresponding to $\beta'$ is larger than $\lambda$, which contradicts the assumption that $\lambda$ is a $(2k+1, 2k+3)$-core partitions into distinct parts which has the largest size. This completes the proof. \qed

 \begin{lemma}\label{lem2.3}
Fix $k\geq 1$. Let $\lambda$ be a $(2k+1, 2k+3)$-core partition into distinct parts such that   $\beta(\lambda)\in S(P_{2k+1, 2k+3})$.    If $\lambda$ is of maximum size, then  there exist some integer   $1\leq i\leq k$  such that   $\beta(\lambda)=\gamma_{k,i}$.
\end{lemma}
\pf  Let $i$ be the maximal integer such that $\beta(\lambda)$ contains the element $2k+2+(i-1)(2k+3)$. Let $j$ be the maximal integer such that $\beta(\lambda)$ contains the element $1+(j-1)(2k+3)$ once $i$ is determined.  By Lemma \ref{22.2}, we have $i< j$. we claim that $\beta(\lambda)$ contains  the element  $1+m(2k+3)$ for all $0\leq m\leq k$. If not, suppose that $\ell$ is the smallest integer such that $\beta(\lambda)$ does not contain the element $1+\ell(2k+3)$. It is easy to check that $\ell> j-1$. By replacing the element $2k+2+(i-1)(2k+3)$ by the element $1+\ell(2k+3)$, we get a new nice order ideal $\beta'$  of $P_{2k+1, 2k+3}$. It is easy to check that $\beta'$ has the same cardinality as $\beta(\lambda)$ and has a larger sum of the elements. By relation (\ref{eq}), the size of the partition corresponding to $\beta'$ is larger than that of $\lambda$. This yields a contradiction with the fact that $\lambda$ has the largest size, which completes the proof of the claim. Combining the claim and Lemmas \ref{22.1} and \ref{22.2}, we have $\beta(\lambda)=\gamma_{k,i}$.
 This completes the proof. \qed

\begin{lemma}\label{lem2.7}
Let $k\geq 2$.  We have $|\mu_{k,i}|\leq |\mu_{k,k}|$ for $1\leq i\leq k$, with the equality holding if and only if $i=k$.
\end{lemma}

\pf By relation (\ref{eq}), the size of $\mu_{k,i}$ is given by
\begin{equation}\label{eq2.11}
\begin{array}{lll}
|\mu_{k,i}|&=& \sum_{h\in \gamma_{k,i}}h-{|\gamma_{k,i}|\choose 2}\\[5pt]
&=& \sum_{h\in \beta_{1,k}}h+\sum_{p=0}^{i-1}(2k+2+p(2k+3))-{|\beta_{1,k}|+i\choose 2}\\[5pt]
&=& \sum_{h\in \beta_{1,k}}h-{|\beta_{1,k}|\choose 2}+ \sum_{p=0}^{i-1}(2k+2+p(2k+3))-{|\beta_{1,k}|+i\choose 2}+{|\beta_{1,k}|\choose 2}\\[5pt]
&=& |\lambda_{1,k}| +(2k+2) {i+1\choose 2}-i|\beta_{1,k}|.
\end{array}
\end{equation}
By the definition of $\beta_{1,k}$, we obtain that
\begin{equation}\label{eq2.12}
|\beta_{1,k}|=
{k+1\choose 2}.
 \end{equation}
Hence, we have
\begin{equation}\label{eq2.13}
|\mu_{k,i}|=|\lambda_{1,k}| +(2k+2) {i+1\choose 2}-i{k+1\choose 2}.
\end{equation}
In particular,   we have
\begin{equation}\label{eq2.14}
|\mu_{k,k }|-|\lambda_{1,k}|= (k+2){k+1\choose 2}={k(k+1)(k+2)\over  2}
\end{equation}
and

 \begin{equation}\label{eq2.15}
|\mu_{k,1}|-|\lambda_{1,k}|=-{k^2-3k-4\over 2}.
\end{equation}

In view of (\ref{eq2.14}) and (\ref{eq2.15}), we deduce that
\begin{equation}\label{eq2.16}
\begin{array}{lll}
|\mu_{k,k }|-|\mu_{k,1}|&=& {k(k+1)(k+2)\over  2}+{k^2-3k-4\over 2} \\
&=&{k^3+4k^2-k-4\over 2}\\
&=& {(k-1)(k+1)(k+4)\over 2}.
\end{array}
\end{equation}

For $k\geq 2$,  it is easy to verify that $ {(k-1)(k+1)(k+4)\over 2}>0$. Hence, for $k\geq 2$,  we have
  \begin{equation}\label{eq2.17}
|\mu_{k,k}|> |\mu_{k,1}|.
\end{equation}

For fixed $k$,  we see that $\mu_{k,i}$ is a quadratic function of $i$ with a positive leading coefficient. Hence the maximum value of $\mu_{k,i}$ is obtained at $i=1$ or $i=k $ when $i$ ranges over $[1,  k ]$. In view of  (\ref{eq2.17}), we have   $|\mu_{k,i}|\leq |\mu_{k,k}|$ for $1\leq j\leq k$, with the equality holding if and only if $i=k$.   This completes the proof.
 \qed

In view of (\ref{eq2.10}) and (\ref{eq2.14}), we have $|\mu_{k,k}|>|\lambda_{2,k}|$ for all $k\geq 2$. One can easily verify that  $|\mu_{1,k}|>|\lambda_{2,k}|$ for $k=1$. Hence, from Lemmas \ref{lem2.4}, \ref{lem2.3},  and \ref{lem2.7}, we deduce the following result.

\begin{theorem}\label{th2.1}
Fix $k\geq 1$. Let $\lambda$ be a $(2k+1, 2k+3)$-core partition into distinct parts which  has the largest size. Then we have $\lambda=\mu_{k,k}$.
\end{theorem}

Now we are ready to derive the largest size of $(2k+1,2k+3)$-core partitions into distinct parts.

\begin{theorem}\label{th2.2}
Let $k\geq 0$.  The largest  size of  $(2k+1,2k+3)$-core partitions into distinct parts is given by
${{k(k+1)(k+2)(5k+11)\over 24}}$.
\end{theorem}
\pf Clearly, the assertion holds for $k=0$. For $k\geq 1$, Let $\lambda$  be a $(2k+1, 2k+3)$-core partition into distinct parts which  has the largest size. By Theorem \ref{th2.1}, we have $\lambda=\mu_{k,k}$.
 In view of (\ref{eq2.14}), we have
 \begin{equation}\label{eq2.18}
 |\lambda|=|\mu_{k,k}|=|\lambda_{1,k}|+{k(k+1)(k+2)\over  2}.
 \end{equation}
By relation (\ref{eq}) and the definition of $\beta_{1,k}$, we have
 \begin{equation}\label{eq2.19}
 \begin{array}{lll}
   |\lambda_{1,k}|&=&\sum_{j=0}^{k-1}(k-j)(2j+1)+ \sum_{j=1}^{k}j(k-j)(2k+3)-{{k+1\choose 2}\choose 2}\\
   \end{array}
    \end{equation}
    By simple computation, from (\ref{eq2.18}) and (\ref{eq2.19}), we can deduce that
    $$
    |\lambda|=|\mu_{k,k}|={{k(k+1)(k+2)(5k+11)\over 24}},
    $$
    as desired. \qed

   Let $s$ be an odd positive integer.  Substituting  $k={s-1\over 2}$ in Theorem \ref{th2.2}, we are led to a proof of Conjecture \ref{con2}.

\noindent{\bf Acknowledgments.} This work was supported by  the National Natural Science Foundation of China (11571320) and  Zhejiang Provincial Natural Science Foundation of China (LY14A010009 and LY15A010008).


\end{document}